\documentclass[11pt,a4paper]{article}
\usepackage{amsmath, amsfonts, amssymb,}
\usepackage{a4wide}
\usepackage{parskip}
\usepackage{enumitem}
\usepackage{xcolor}
\usepackage{indentfirst}
\usepackage{pstricks}
\usepackage{latexsym}
\usepackage{amsthm}
\usepackage{cite}

\newcommand{\e}{\epsilon}

\newcommand{\al} {\alpha}
\newcommand{\ba} {\beta}
\newcommand{\de} {\delta}
\newcommand{\ga} {\gamma}

\newcommand{\De} {\Delta}
\newcommand{\la} {\lambda}

\newcommand{\noi} {\noindent}

\newcommand{\ds} {\displaystyle}

\newcommand{\RR}{{\mathbb R}}
\newcommand{\te} {\theta}

\newtheorem{thm}{Theorem}[section]
\newtheorem{rem}{Remark}[section]
\newtheorem{lem}{Lemma}[section]
\newtheorem{prop}{Proposition}[section]

\numberwithin{equation}{section}

\begin{document}
	\setlength{\abovedisplayskip}{3pt}
	\setlength{\belowdisplayskip}{3pt}
	\date{}
	{\vspace{0.01in}
		\title{Existence of positive solutions for a class of quasilinear Schr\"{o}dinger equations with critical Choquard nonlinearity}

	\author{ { Sushmita Rawat \footnote{e-mail: {\tt sushmita.rawat1994@gmail.com}} and 
			K. Sreenadh \footnote{e-mail: {\tt sreenadh@maths.iitd.ac.in}}} \\ 
		Department of Mathematics, Indian Institute of Technology Delhi,\\ Hauz Khas, New Delhi 110016, India.}

	\maketitle
	\begin{abstract}
		\noi This article is concerned with the existence of positive weak solutions for the following quasilinear Schr\"odinger Choquard equation:
		\begin{equation*}
			\begin{array}{cc}
				\ds  -div(g^2(u)\nabla u) + g(u)g'(u)|\nabla u|^2 + a(x) u = k(x, u) \;\text{in} \; \RR^N,
			\end{array}
		\end{equation*}
		where $N \geq 3$, $\ds k(x,u) := h(x,u) + (I_{\vartheta}*|u|^{\al\cdot2^*_\mu})|u|^{\al\cdot2^*_\mu-2}u$, $g : \RR \to \RR^+$ is a differentiable even function with $g(0) = 1$ and $g'(t) \geq 0$ for all $t \geq 0$
		; $h\in C( \RR^N\times\RR, \RR)$ and the potential $a \in C( \RR^N, \RR)$. We establish the existence of positive solution using the change of variable and variational methods, under appropriate assumptions on $g$, $h$ and $a$.

		\medskip

		
		\medskip
		

	\end{abstract}
	\newpage
	\section{Introduction}
	In this work, we study the existence of positive solution for the quasilinear Schr\"odinger equation with Hardy-Littlewood-Sobolev critical exponent. Precisely we consider the problem
	\begin{equation*}    (P)\; \left\{  \begin{array}{cc}
			-div(g^2(u)\nabla u) + g(u)g'(u)|\nabla u|^2 + a(x) u = k(x, u) \;\text{in} \; \RR^N,
		\end{array} \right.
	\end{equation*}
 where we assume $\ds k(x,u) = h(x,u) + (I_{\vartheta}*|u|^{\al\cdot2^*_\mu})|u|^{\al\cdot2^*_\mu-2}u$, \; $x \in \RR^N$; $0 <\vartheta < N$; $I_{\vartheta} (x) = \dfrac{A_{\vartheta}}{|x|^{\mu}}$ is the Riesz potential and $A_{\vartheta}$ is an appropriate constant; $0 < \mu = N-\vartheta$. Here $2^{*}_{\mu } = \frac{2N-\mu}{N-2}$ is the critical exponent in the sense of Hardy-Littlewood-Sobolev inequality and we assume $\mu < \min\{N, 4\}$. 
 
Our motivation to study the problem $(P)$ mainly comes from the fact that the solutions of $(P)$ are related to the existence of solitary wave
 solutions for quasilinear Schr\"odinger equations
 \begin{equation}\label{Deq1.1}
 	i\partial_tz = -\De z + W(x)z - k(x, z) - \De(l(|z|^2))l'(|z|^2)z,\; x \in\RR^N,
 \end{equation}
 where $z : \RR \times \RR^N \to \mathbb{C}$, $W : \RR^N \to \RR$ is a given potential, $l$ and $k(x, z)$ are suitable real functions. 
 Putting $z(t, x) = e^{(-iEt)}u(x)$ in \eqref{Deq1.1}, where $E \in \RR$ and $u > 0$ is a real function, we attain an analogous elliptic equation
 \begin{equation}\label{Deq1.2}
 	-\De u + a(x)u -\De(l(u^2))l'(u^2)u = k(x,u),\; x \in \RR^N.
 \end{equation}
For the case when $\ds g^2(u) = 1 + \frac{[(l(u^2))']^2}{2},$ then Problem $(P)$ and \eqref{Deq1.2} are equivalent.

 The quasilinear equations of the form \eqref{Deq1.1},
 have been accepted as models of several physical phenomena corresponding to various types of nonlinear terms $l$. For instance, Kurihura in \cite{kurihura}, studied the non-linear time evolution of the condensate wave function in the superfluid films, for $l(s) = s$. 
 The equation \eqref{Deq1.2} with $l(s) = s$, has been researched extensively with a range of conditions on the potential $a(x)$ and 
 with subcritical nonlinearity $k(x, u)= |u|^{q-2}u$, $4 \leq q < 2\cdot2^*$,  as observed in \cite{liu3}, the number $2\cdot2^*$ behaves like a critical exponent, where $2^* = \frac{2N}{N-2}$. Authors in \cite{liu2,jeanjean}, introduced a new formulation of the problem by using a change of variables which proved to be an efficient way to examine the existence and multiplicity of solutions. See for example \cite{liu1, liu2, liu3, poppenberg, jeanjean, silva1, rad_zhang} and the reference therein.
 Whereas for the critical case,
 existence of solutions was first studied by Moameni in \cite{moameni} when the potential function $a(x)$ satisfies some geometry conditions. 
 For more results on the critical case we cite \cite{miyagaki, silva2, wang, lins, liu4} and references therein.
 Furthermore the quasilinear equations \eqref{Deq1.1}, with $l(s) = (1 + s)^{\frac{1}{2}}$, 
models the self-channeling of a high-power ultra short laser in matter. We refer to \cite{J_yang, de.hayashi.saut, huang} and so on.\\
In recent years, quasilinear equation \eqref{Deq1.2}, with general function $l$ has received much attention.
Shen and Wang \cite{shen_wang}, introduced a new variable replacement and obtained the existence of positive solutions
for $(P)$ when $k(x, u)$ is superlinear and subcritical. Later, under some suitable assumptions on g, k, and a, the result was extended by Deng et al. in \cite{deng2, deng3}, and they established the existence of positive solutions when $k(x, u)$ is critical. There have been several studies about the existence of nontrivial solutions one can further refer to \cite{ shi_chen2, li_Wu, quanqing_xian} and so on.

On the other hand, for $g(t) =1$, we deduce that $(P)$ corresponds to the following semilinear problem
\begin{equation}\label{Deq1.5}
	-\Delta u + a(x) u = (|x|^{-\mu}* F(u))f(u) \text{ in } \RR^N
\end{equation}
where $f \in C(\mathbb{R}, \mathbb{R})$ satisfies some growth condition,   $F$ is  the anti-derivative of $f$ and $a$ is the vanishing potential.
A comprehensive study has been done on the existence and uniqueness results of the above mentioned nonlocal elliptic equations due to their vast applications in physical models. One of the first applications of Choquard equations was given by Pekar in the framework of quantum theory \cite{pekar} and Lieb \cite{choqlieb} used it in the approximation of Hartree-Fock theory.
 Moroz and Schaftingen in \cite{Moroz3}, studied \eqref{Deq1.5}. For a detailed state of the art research, readers can refer \cite{Moroz4, M.yang} and references therein.\\
The equation of type \eqref{Deq1.5} is usually called the nonlinear Schr\"odinger-Newton equation. If $u$ solves \eqref{Deq1.5} then the function $z$ defined by $z(t,x) = e^{it}u(x)$ is a solitary wave of the focusing time-dependent Hartree equation
\begin{equation*}
	i\partial_tz = -\De z + W(x)z - \left(I_\al* |z|^p\right)|z|^{p-2}z \;\text{ in } \RR^+ \times \RR^N.
\end{equation*}

On the contrary, results about the combination of problem \eqref{Deq1.2} and problem \eqref{Deq1.5} are relatively few. Yang et al. in \cite{x_yang} considered the quasilinear equation \eqref{Deq1.2} with $l(s) = s$ and $k(x,u) = (|x|^{-\mu} * |u|^p)|u|^{p-2}u$ in $\RR^N$,
where $N \geq 3$, $\mu \in (0, (N + 2)/2)$, $p \in (2, (4N - 4\mu)/(N - 2))$. By using the perturbation method, they obtained the existence of positive solutions, negative solutions, and high energy
solutions. We also cite \cite{s.chen_x.wu, x_yang.x_tang, s.liang_y.song, s.liang_b.zhang} and the references therein.
However, for general function $l$, the existence of solutions to the generalized quasilinear Schr\"odinger equation involving critical Choquard-type term is still
open, to the best of our knowledge.

 \noindent Motivated by the works described above, in the present paper we are concerned with the existence of positive solutions to our problem $(P)$. Due to the appearance of the critical term and the domain $\RR^N$, the compactness of the embedding is lost. Moreover, the characteristics of the term $\ds F^*(s) = \frac{1}{\al 2^*}(G^{-1}(s))^{\al 2^*} - \frac{\al^{2^*-1}}{2^*\ba^{2^*}}|s|^{2^*}$, obtained for the critical case as in Deng et al. \cite[Lemma 2.2]{deng3}, can not be extended for the Choquard case. To conquer these difficulties, we first focus on the limiting case and employ the concentration-compactness principle to overcome the problem with the help of assumption $(g_1)(b)$. Following, we prove some delicate estimates concerning the critical Choquard term to prove that the mini-max value obtained along the Mountain-pass theorem is smaller than the admissible threshold for the Palais-Smale condition. Thus we obtain a positive radial solution, which is also a ground-state solution. Then using these estimates and the critical level of the functional $J^\infty$, we prove the existence of a positive solution for our problem $(P)$.
 
Throughout the paper, we assume $h: \RR^N\times \RR \to \RR$ is a continuous function
  and denote $\ds H(x,u) = \int\limits_{0}^{u} h(x,t)dt$. We observe that the energy functional 
 \begin{equation*}
 	I(u) = \frac{1}{2}\int\limits_{\RR^N} g^2(u)|\nabla u|^2 + a(x)u^2dx -\int\limits_{\RR^N}H(x,u)dx -
 	\frac{1}{2\al\cdot{2^*_{\mu } }}\iint\limits_{\RR^{2N}} \frac{|u(y)|^{\al\cdot2^{*}_{\mu }}|u(x)|^{\al\cdot2^{*}_{\mu }}}{|x-y|^ \mu}dxdy
 \end{equation*}
 associated with $(P)$ might not be well defined in $H^1(\RR^N)$. 
  In the spirit of the argument developed by Shen and Wang in \cite{shen_wang}, we make the change of variables, 
 $$v= G(u) = \int\limits_{0}^{u} g(t)dt. $$ 
 Then we have
 \begin{equation*}
 	\begin{aligned}
 		J(v) = & \frac{1}{2}\int\limits_{\RR^N} |\nabla v|^2 dx +\frac{1}{2}\int\limits_{\RR^N}a(x)(G^{-1}(v))^2dx -\int\limits_{\RR^N}H(x,G^{-1}(v))dx \\
 		&-\frac{1}{2\al\cdot{2^*_{\mu } }}\iint\limits_{\RR^{2N}} \frac{|G^{-1}(v(y))|^{\al\cdot2^{*}_{\mu }}|G^{-1}(v(x))|^{\al\cdot2^{*}_{\mu }}}{|x-y|^ \mu}dxdy,
 	\end{aligned}
 \end{equation*}
 which is well defined on the usual Sobolev space $H^1(\RR^N )$, under some suitable assumptions on the functions $a$, $g$, and $h$.\\
 Further, taking into account that we are looking for a positive solution, we rewrite the functional $J$ as
 \begin{equation*}
 	\begin{aligned}
 		\ds J(v) = & \frac{1}{2}\int\limits_{\RR^N} |\nabla v|^2 dx +\frac{1}{2}\int\limits_{\RR^N}a(x)(G^{-1}(v))^2dx -\int\limits_{\RR^N}H(x,G^{-1}(v))dx \\
 		&-\frac{1}{2\al\cdot{2^*_{\mu } }}\iint\limits_{\RR^{2N}} \frac{|G^{-1}(v^+(y))|^{\al\cdot2^{*}_{\mu }}|G^{-1}(v^+(x))|^{\al\cdot2^{*}_{\mu }}}{|x-y|^ \mu}dxdy.
 	\end{aligned}
 \end{equation*}
 The functional $J \in C^1$, indeed for $\phi \in C^{\infty}_0(\RR^N)$
 \begin{equation*}
 	\begin{aligned}
 		\ds \langle J'(v), \phi\rangle= & \int\limits_{\RR^{N}}\nabla v\cdot\nabla\phi + a(x)\frac{G^{-1}(v)}{g(G^{-1}(v))}\phi -\frac{h(x,G^{-1}(v))}{g(G^{-1}(v))}\phi dx\\
 		& - \iint\limits_{\RR^{2N}}\frac{|G^{-1}(v^+(y))|^{\al\cdot2^{*}_{\mu }}|G^{-1}(v^+(x))|^{\al\cdot2^{*}_{\mu}-2}G^{-1}(v^+(x))\phi(x)}{|x-y|^ \mu \,g(G^{-1}(v^+(x)))}dxdy.
 	\end{aligned}
 \end{equation*}
 Thus our problem reduces to examining the semilinear equation
 \begin{equation}\label{Deq1.6}
 		-\De v + a(x) \frac{G^{-1}(v)}{g(G^{-1}(v))} = \frac{h(x, G^{-1}(v))}{g(G^{-1}(v))} + \left(I_{\vartheta}*|G^{-1}(v^+)|^{\al\cdot2^*_\mu}\right)|G^{-1}(v^+)|^{\al\cdot2^*_\mu-2}\frac{G^{-1}(v^+)}{g(G^{-1}(v^+))}, 
 \end{equation}
 
Before stating our main results we make some assumptions on the functions $ g, h, \;\text{and}\; a$.\\
We assume $g$ satisfies the following conditions:
\begin{enumerate}
	\item[$\mathbf{(g_0)}$] $g\in C'(\RR, \RR^+)$, $g(t)$ is even, $g(0) = 1$ and $g'(t) \geq 0$ for all $t\geq 0$;
	\item[$\mathbf{(g_1)}$] There exist constants $\al \geq 1$, $\ba >0$ and $\gamma \in (-\infty, \al)$, such that for $t\to +\infty$,
	${g(t) = \ba t^{\al-1} + O(t^{\ga -1})}$. Further
	\begin{enumerate}
		\item $(\al -1)g(t) \geq g'(t)t \; \text{for all} \; t\geq 0$;
		\item $g'(t) \geq \ba^2(\al -1) \frac{t^{2\al-3}}{g(t)}.$
	\end{enumerate} 
	\end{enumerate}
\begin{rem} Functions $g(t) = \sqrt{1 + \frac{(q+1)^2}{2}t^{2q}}$, for $q >\frac{1}{2}$ naturally satisfies conditions $(g_0), (g_1)$, for $\al = q + 1$ and $\ba = \frac{q+1}{\sqrt{2}}$. In particular for $q =1$ we have
$g(t) = \sqrt{1+2t^2}$. 
Also one can consider functions of the type $g(t) = (1+ t^2)^{\frac{\al -1}{2}}$.
\end{rem}
Next we assume $h$ and $a$ satisfy the following conditions, where $2^*= \dfrac{2N}{N-2}$,
\begin{enumerate}
	\item [$\mathbf{(h_0)}$] $ h(x,t) \geq 0$, differentiable in $t \in [0, \infty) $ for all $x \in \RR^N$ and continuous in $x \in \RR^N$ for all $t \in [0, \infty) $. For all $t \leq 0 $ and $x \in \RR^N$, $h(x,t) \equiv 0$;
	\item[$\mathbf{(h_1)}$] $\lim\limits_{t \to +\infty}\dfrac{h(x,t)}{|t|^{\al\cdot 2^* -1}} =0$ and $\lim\limits_{t \to 0}\dfrac{h(x,t)}{t} = 0$ uniformly in $x \in \RR^N$;
	\item[$\mathbf{(h_2)}$] There exists $\tilde{\mu} \in (2, 2^*)$ such that for any $t >0$
	\begin{equation*}
		\frac{\partial h(x,t)}{\partial t}t \geq (\al \tilde{\mu} -1)h(x,t);
	\end{equation*}
	\item[$\mathbf{(h_3)}$] $\lim\limits_{|x| \to +\infty}h(x,t) = \overline{h}(t) $ uniformly on any compact subset of $[0, \infty)$  and there exists a constant $\nu > 2$, such that for any $\e > 0$ we can find $C_\e > 0$ satisfying
	$$h(x,t) - \overline{h}(t) \geq -e^{-\nu|x|}(\e t^\al + C_\e t^{\al p})t^{\al-1} \;\text{for all}\, x\in \RR^N, t \geq 0, p \in (1, 2^*-1 );$$
\end{enumerate}
\begin{rem}
	Functions of the type,
	$h(x,t)$ = $\begin{cases}(1- e^{-\nu|x|})t^{\al q+\al-1} & \text{for}\; t \geq 0,\; x \in \RR^N,\\
		0 &\;\text{for}\; t < 0,\; x \in \RR^N,
	\end{cases}$
	satisfy our assumptions $(h_0)-(h_3)$, where $\nu > 2$, $1 <q <p$ and $p\in (1, 2^*-1 )$.
\end{rem}
\begin{enumerate}
	\item[$\mathbf{(a_0)}$] $0 \leq a(x) \in C(\RR^N ), \lim\limits_{|x| \to +\infty}a(x) = 1$, $a(0) < 1$, and $0 \leq 1 - a(x) \leq k e^{-\nu|x|}$ for some positive
	constant $k$.
\end{enumerate}

\begin{rem}\label{Drem1.1} 
	Using the fact that $\dfrac{g(s)}{|s|^{\al-1}}$ is a decreasing function, we claim for all $s \in \RR$
$$g(s) \geq \ba |s|^{\al -1}.$$ 
\end{rem}

\begin{rem}
	Deng et al. in \cite{deng3}, found that $\al\cdot2^*$ is the critical exponent for problem
	$(P)$ if $\lim \limits_{s\to + \infty}\frac{g(s)}{s^{\al-1}}:= \ba > 0$, for some $\al \geq 1$. Moreover, they proved that $(P)$ with $k(x, u) \equiv 0$ in $\RR^N\times\RR$ has no positive solution, if $p \geq \al\cdot2^*$ and
	$x \cdot\nabla a(x) \geq 0$ in $\RR^N$.
	
\end{rem}

With this introduction, we state our main result:
\begin{thm}\label{Dthm1.1}
	Suppose  $(g_0)-(g_1)$, $(h_0)-(h_3)$ and $(a_0)$ hold. Then for the problem \eqref{Deq1.6} which is equivalent to $(P)$, there exists at least one positive solution provided $N$ and $\tilde{\mu}$ satisfies: 
	\begin{enumerate}
		\item[(1)] $N \geq \max\left\lbrace 2 + \frac{4\al}{\al(\tilde{\mu}-1)-\ga^+},\, 6\right\rbrace$ and $ \tilde{\mu} > 2$,
		\item [(2)] $N = 5;\, \tilde{\mu} > \frac{7}{3} +\frac{\ga^+}{\al},$
		\item [(3)] $N = 4;\, \tilde{\mu} > 3 +\frac{\ga^+}{\al},$
		\item [(4)] $N = 3;\, \tilde{\mu} > 5 +\frac{\ga^+}{\al},$
	\end{enumerate} where $\ga^+ = \max \{\ga, 0\}$.
\end{thm}
We prove this with the help of limiting equation at infinity
\begin{equation*}  (P^\infty)\; \left\{  \begin{array}{cc}\ds
		-\De v + \frac{G^{-1}(v)}{g(G^{-1}(v))} =  \frac{\overline{h}(G^{-1}(v))}{g(G^{-1}(v))}  + \left( I_{\vartheta}*|G^{-1}(v^+)|^{\al\cdot2^*_\mu}\right) \frac{|G^{-1}(v^+)|^{\al\cdot2^*_\mu-2}G^{-1}(v^+)}{g(G^{-1}(v^+))}, 
	\end{array} \right.
\end{equation*}
and the associated functional $J^{\infty}$ is defined as
\begin{equation*}
	\begin{aligned}
		\ds J^{\infty}(v) = & \frac{1}{2}\int\limits_{\RR^N} |\nabla v|^2  +(G^{-1}(v))^2dx -\int\limits_{\RR^N}\overline{H}(G^{-1}(v))dx \\
		&-\frac{1}{2\al\cdot{2^*_{\mu } }}\iint\limits_{\RR^{2N}} \frac{|G^{-1}(v^+(y))|^{\al\cdot2^{*}_{\mu }}|G^{-1}(v^+(x))|^{\al\cdot2^{*}_{\mu }}}{|x-y|^ \mu}dxdy.
	\end{aligned}
\end{equation*}
 
Throughout the paper we will use the following notations 
$u^+ = \max\{u, 0\}$, 
\begin{equation*}
		\|u\|_{0}^{2\cdot2^{*}_{\mu}} := \iint\limits_{\RR^{2N}} \frac{|u(x)|^{2^{*}_{\mu}}|u(y)|^{2^{*}_{\mu}}}{|x-y|^ \mu}\,dxdy\quad \text{and} \quad
	\|u\|_{0,\al}^{2\al\cdot2^{*}_{\mu}} := \iint\limits_{\RR^{2N}} \frac{|u(x)|^{\al\cdot2^{*}_{\mu}}|u(y)|^{\al\cdot2^{*}_{\mu}}}{|x-y|^ \mu}\,dxdy.
\end{equation*}
The letters $C$, $C_i$ will denote various positive constants
whose exact values may change from line to line but are inessential to the analysis of the problem.

Rest of the paper is organized as follows. In Section 2, we present some preliminary knowledge and set up the functional. In Section 3, we give some technical lemmas which will help in proving our Theorem \ref{Dthm1.1}. In Section 4, we work on the limiting case and prove it has a positive radial solution. Section 5 is devoted to the proof of
Theorem \ref{Dthm1.1}.

\section{Preliminaries}
The functional space associated to this problem is $H^1(\RR^N)$ with the corresponding norm,
\begin{equation*}
	\|u\|^2 = \ds \int\limits_{\RR^{N}} |\nabla u|^2\,dx + \int\limits_{\RR^{N}} | u|^2\,dx.
\end{equation*}
	Next we recall the famous Hardy-Littlewood-Sobolev inequality,
\begin{prop}
	\textbf{(Hardy-Littlewood-Sobolev inequality)}: Let $t$, $r > 1$ and $0 < \mu < N$ with $\frac{1}{t} + \frac{\mu}{N} + \frac{1}{r} = 2$, $f \in L^t(\mathbb{R}^N)$ and $h \in L^r(\mathbb{R}^N)$. Then there exists a sharp constant $C(t, r, \mu, N)$ independent of $f$, $h$ such that
	\begin{equation*}
		\iint\limits_{\mathbb{R}^{2N}} \dfrac{f(x)h(y)}{|x-y|^ \mu}\,dxdy \leq C(t, r, \mu, N)\|f\|_{L^t(\mathbb{R}^N)}\|h\|_{L^r(\mathbb{R}^N)}.
	\end{equation*}
	\end{prop}

From the classical Sobolev embedding theorem, we have that $H^1(\RR^N)$ is continuously embedded into $L^{p}(\RR^N)$, for all $ p \in [2, 2^{*}]$. 
The best constant for the embedding $H^1(\RR^N)$ into $L^{2^{*}}(\mathbb{R}^N)$ is 
\begin{equation}\label{Deq2.1}
	S =\inf\limits_{u\in H^1(\RR^N)\backslash \{0\}}\left\lbrace\int\limits_{\mathbb{R}^{N}} |\nabla u|^2\,dx : \int\limits_{\mathbb{R^N}}|u|^{2^{*}} = 1 \right\rbrace.
\end{equation}
Consequently, we define
\begin{equation}\label{Deq2.2}
	S^H = \inf\limits_{u\in H^1(\RR^N)\backslash \{0\}}\left\lbrace\int\limits_{\mathbb{R}^{N}} |\nabla u|^2\,dx:  \iint\limits_{\mathbb{R}^{2N}} \frac{|u(x)|^{2^{*}_{\mu}}|u(y)|^{2^{*}_{\mu}}}{|x-y|^ \mu}\,dxdy = 1\right\rbrace. 
\end{equation}
We shall summarize the notion and notations of the function where the infimum of  \eqref{Deq2.1} and \eqref{Deq2.2} exists.
\begin{lem} \cite[Lemma 1.2]{M.yang}
	The constant $S^H$ is achieved by u if and only if the form of $u$ is
	$C\left( \frac{t}{t^2 + |x-x_0|^2}\right) ^\frac{N-2}{2}$, $x \in \mathbb{R}^N$, 
	for some $x_0 \in \mathbb{R}^N, C \text{and}\; t > 0.$ Moreover,
	$S^H = \dfrac{S}{{C(N, \mu)}^\frac{1}{2_{\mu}^{*}}}$.
	\qed
\end{lem}

Consider the family of functions ${U_\epsilon}$, where $U_\epsilon$ is defined as
\begin{equation} \label{minmimizer}
	U_\epsilon = \epsilon^{-\frac{N-2}{2}}u^{*}\left( \frac{x}{\epsilon}\right), \; x\in \mathbb{R}^{N}, \epsilon > 0,
\end{equation}
\begin{equation*}
	u^{*}(x) = \overline{u}\left( \frac{x}{{S}^\frac{1}{2}}\right) , \; \overline{u}(x) = \frac{\tilde{u}(x)}{\|\tilde u\|_{L^{2^{*}}(\mathbb{R}^N)}} \; \text{and}\; \tilde{u}(x) = \alpha(\beta^2 + |x|^2)^{-\frac{N-2}{2}},
\end{equation*}
with $\alpha > 0$ and $\beta > 0$ are fixed constants.
Then for each $\epsilon > 0, \; U_\epsilon $ satisfies
\begin{equation*}
	-\Delta u = |u|^{2^{*}-2}u \quad in \; \mathbb{R}^N,
\end{equation*}
and the equality,
\begin{equation*}
	\int\limits_{\mathbb{R}^{N}} |\nabla U_\epsilon|^2\,dx = \int\limits_{\mathbb{R^N}}|U_\epsilon|^{2^{*}} = {S}^\frac{N}{2 }. 
\end{equation*}
Let $\eta \in C_0^{\infty}(\mathbb{R}^N,[0,1])$ be a radial cut-off function such that $\eta(|x|)$ =	$\begin{cases}
		1 & \quad B_{\rho_\e},\\
		0 & \quad \mathbb{R}^N\backslash B_{2\rho_\e},
	\end{cases}$
where $\rho_\e = \e^\tau$, $\tau \in ( \frac{1}{2}, 1)$
and for each $\epsilon > 0$, set 
\begin{equation}\label{Deq2.4}
	u_\epsilon(x) = \eta(x)U_\epsilon(x)\quad \text{for} \; x\in \mathbb{R}^N,
\end{equation}
where $U_\epsilon$ is as defined in \eqref{minmimizer}.

\begin{prop}\label{Dprop2.2} \cite[Theorem 3.3]{s.bae}, \cite[Proposition 2.8]{Tuhina} Let $N \geq 3$. Then 
	\begin{enumerate}
		\item[(1)]
			$\|\nabla u_\epsilon\|_{L^2}^2 = C^\frac{N}{2\cdot2^*_{\mu}}(S^H)^\frac{N}{2} + O(\epsilon^{N-2}), \; \|u_\epsilon\|^{2^{*}}_{L^{2^{*}}} = C^\frac{N}{2\cdot2^*_{\mu}}(S^H)^{\frac{N}{2}}+ O(\epsilon^N),$
	\item[(2)]$ 
		\|u_\epsilon\|^2_0\geq \left[C(N, \mu)^{\frac{N}{2}}(S^H)^{\frac{2N-\mu}{2}}-O(\epsilon^N)\right]^{\frac{1}{2^*_{\mu}}},$
\item[(3)] $\ds
	\int\limits_{\RR^{N}} |u_\epsilon|^t = \begin{cases}
		C\e^{N-\frac{t(N-2)}{2}}+ O(\e^{\frac{t(N-2)}{2}}) & t > \frac{N}{N-2},\\
		C\e^{\frac{t(N-2)}{2}}|\ln(\e)|+ O(\e^{\frac{t(N-2)}{2}}) & t = \frac{N}{N-2},\\
		O(\e^{\frac{t(N-2)}{2}}) & t < \frac{N}{N-2}\\
	\end{cases}$
\end{enumerate}
as $\epsilon \to 0$, for some positive constant $C$.
	
\end{prop}

\section{Technical lemmas}
In this section, we gather a couple of critical results, which will be used to prove our main results. First, we present some properties of the functions $g$ and $h$:
\begin{lem}\label{Dlem3.1}
	The functions $h(x,t)$, $g(t)$ and $H(x,t) = \int_{0}^{t}h(x,\tau) d\tau$, $G(t) = \int_{0}^{t}g(\tau)d\tau$ satisfy the
	following properties under the assumptions $(g_0), (g_1)$, and $(h_0)-(h_2)$
\begin{enumerate}
	\item[(1)] $G(t)$ and $G^{-1}(s)$ are odd.
	\item[(2)] For all $t \geq 0$, $s \geq 0$,
	$G(t) \leq g(t)t$ and $G^{-1}(s) \leq \dfrac{s}{g(0)}$.
	\item[(3)] For all $s \geq 0$, $\dfrac{G^{-1}(s)}{s}$ is non-increasing and\\
	$\lim\limits_{s \to 0}\dfrac{G^{-1}(s)}{s} = \dfrac{1}{g(0)}$,\, $\lim\limits_{s \to \infty}\dfrac{G^{-1}(s)}{s} = \begin{cases}	\dfrac{1}{g(\infty)} & \text{if}\; g \;\text{is bounded};\\
	0 & \text{if} \;g\; \text{is unbounded}.
\end{cases}$
\item[(4)] For any $x \in \RR^N$ , $t > 0$, there holds\\
$\al G(t) \geq g(t)t$,\, $\al \tilde{\mu} H(x,t) \leq h(x,t)t$; \\
$h(x,t)G(t) \geq \tilde{\mu}g(t)H(x,t)$,\, $(\tilde{\mu}- 1)h(x,t) \leq G(t)\left(\dfrac{h(x,t)}{g(t)}\right)'_t$;
\item[(5)] There exist constants $C > 0$, $M > 0$ such that $H(x,t) \geq C(G(t))^{\tilde{\mu}}$ for all $x \in \RR^N$ , $t \geq M$.
\item[(6)] $t^\al \leq \dfrac{\al}{\ba}G(t)$ and there exist some constants $C > 0$ and $M > 0$ such that
\begin{equation*}
	0 \geq \left( \dfrac{t^\al}{G(t)}\right)^{2^*} - \left( \dfrac{\al}{\ba}\right)^{2^*} \geq -C G(t)^{-\de}
\end{equation*}  
for $t \geq M$, where $\de =1-\frac{\ga^+}{\al}$ and $\ga^+ = \max \{\ga, 0\}$.
\end{enumerate}
\proof The proof is immediate, one can refer to \cite[Lemma 2.1]{deng3}.
\qed
\end{lem}
Denote
	\begin{align*}
		f(x,s) &= a(x)s - a(x)\dfrac{G^{-1}(s)}{g(G^{-1}(s))}+ \dfrac{h(x,G^{-1}(s))}{g(G^{-1}(s))},\\ 
		F(x,s) &= \int_{0}^{s}f(x,\tau)d\tau = \frac{1}{2}a(x)[s^2 - (G^{-1}(s))^2] + H(x,G^{-1}(s)).
	\end{align*}
\begin{lem}\label{Dlem3.2}
	The functions $f(x,s)$, $F(x,s)$ satisfy the following properties
	under the assumptions $(g_0)-(g_1)$, $(h_0)-(h_3)$ and $(a_0)$.
	\begin{enumerate}
		\item[(1)] $f(x,s) \geq 0$ for all $x \in \RR^N$, $s \geq 0$.
		\item[(2)] $\lim\limits_{s \to 0^+}\dfrac{f(x,s)}{s}= 0$, $\lim\limits_{s \to 0^+}\dfrac{F(x,s)}{s^2}= 0$ uniformly in $x \in \RR^N$.
		\item[(3)] $\lim\limits_{s \to +\infty}\dfrac{f(x,s)}{s^{2^*-1}} = 0$, $\lim\limits_{s \to +\infty}\dfrac{F(x,s)}{s^{2^*}} = 0$ uniformly in $x \in \RR^N$.
		\item[(4)] $f(x,s)s \geq 2F(x,s)$ for all $x \in \RR^N$, $s \geq 0$.
		\item[(5)] $\lim\limits_{|x| \to +\infty}f(x,s) = \overline{f}(s)$ exist and
		$\overline{f}(s) = s - \dfrac{G^{-1}(s)}{g(G^{-1}(s))}+ \dfrac{\overline{h}(G^{-1}(s))}{g(G^{-1}(s))}$.
		Moreover, for some positive constant $M_1$
		\begin{equation*}
			f(x,s) - \overline{f}(s) \geq -2Ce^{-\zeta|x|}(\e s + C_\e s^p), \;\text{for all}\; x \in \RR^N, s \geq M_1.
		\end{equation*}
	\end{enumerate}
\proof The proof is direct, one can refer to \cite[Lemma 2.2]{deng3}.
\qed
\end{lem}

\begin{lem}\label{Dlem3.5}
	\cite[ Lemma 2.1 ]{zhu_cao} Let $\{\rho_n\}_{n\geq1}$ be a sequence in $L^1(\RR^N )$ satisfying $\rho_n \geq 0$ on $\RR^N$ , $\lim\limits_{n \to \infty}\int\limits_{\RR^{N}}\rho_n\,dx = \la >0,$
where $\la > 0$ is fixed. Then there exists a sub-sequence $\{\rho_{n_k} \}$ satisfying one of the following two possibilities:
	\begin{itemize}
		\item (Vanishing):$\ds\lim\limits_{k \to +\infty}\sup\limits_{y \in \RR^N}\int\limits_{B_R(y)}\rho_{n_k} (x)dx = 0$, for all $0 <R< +\infty$. 
		\item (Nonvanishing): There exist $\al > 0$, $0 <R< +\infty$ and $\{y_k\} \subset \RR^N$ such that\\ $\ds\lim\limits_{k \to +\infty}\int\limits_{B_R(y_k)}\rho_{n_k} (x)dx \geq \al > 0$.
	\end{itemize}
\end{lem}
\begin{lem}\label{Dlem3.6}
	\cite[ Lemma 2.3 ]{zhu_cao} Let $1 < p \leq \infty$, $1 \leq q < \infty$, with $q \neq \dfrac{Np}{N-p}$ if $p<N$. Assume that $v_n$ is bounded in $L^q (\RR^N )$, $|\nabla v_n|$ is bounded in $L^p(\RR^N )$ and
	$\ds\sup \limits_{y\in \RR^N}\int\limits_{B_R(y)}|v_n|^q\;dx \to 0$, for some $R > 0$ as $n \to \infty$.
	Then $v_n \to 0$ in $L^\al(\RR^N )$, for $\al \in (q, \frac{Np}{N-p}). $
	\end{lem}
Next, we look at the Mountain-pass geometry of the functional $J$
\begin{lem}\label{Dlem3.7} The functional $J$ satisfies the following conditions:
	\begin{itemize}
		\item [(i)]There exist $\mathfrak{a},\rho >0$ such that $J(v) > \mathfrak{a}$ for $\|v\| =\rho$.
		\item [(ii)]$J(0) = 0$ and there exists $e \in H^1(\RR^N)$ with $\|e\| > \rho$ and $J(e) < 0$.
	\end{itemize}
	\proof (i). Using Lemma \ref{Dlem3.1}(6), Lemma \ref{Dlem3.2}(2) (3) and equivalence of the norms, we infer that for every $\e>0$, there exists $C_\e>0$ such that
	\begin{equation*}
\begin{aligned}
			J(v) &= \frac{1}{2}\int\limits_{\RR^N} |\nabla v|^2 + a(x)v^2dx  -\int\limits_{\RR^N}F(x,v)dx -\frac{1}{2\al\cdot{2^*_{\mu } }}\|G^{-1}(v^+)\|_{0,\al}^{2\al\cdot2^{*}_{\mu }}\\
			 & \geq \left( \frac{C}{2}-\e\right) \int\limits_{\RR^N} |\nabla v|^2 + v^2dx - C_\e \int\limits_{\RR^{N}}|v|^{2^*}\,dx-\frac{1}{2\cdot2^*_{\mu}}\frac{\al^{2\cdot2^*_\mu-1}}{\ba^{2\cdot2^*_{\mu}}}\|v^+\|_0^{2\cdot2^*_{\mu}},
			\end{aligned}
	\end{equation*}
for some positive constant $C>0$.
	Further using Sobolev embedding \eqref{Deq2.1} and \eqref{Deq2.2}, we get
	\begin{equation*}
		J(v) \geq \left( \frac{C}{2}-\e\right)  \|v\|^2 - S^{\frac{-2^*}{2}} C_\e\|v\|^{2^*}-(S^H)^{-2^*_\mu}\frac{1}{2\cdot2^*_{\mu}}\frac{\al^{2\cdot2^*_\mu-1}}{\ba^{2\cdot2^*_{\mu}}}\|v\|^{2\cdot2^*_{\mu}}.
	\end{equation*}
It is easy to see that for $\|v\|$ small enough, the functional $J$ satisfies the desired geometry.
\proof (ii). Let $\phi \in C^\infty_0(\RR^N)$ with $0 \leq \phi \leq 1$ and supp $\phi = \overline{B_1}$. Then from Lemma \ref{Dlem3.1}(3), we infer that $G^{-1}(t\phi) \geq G^{-1}(t)\phi$, for all $t \geq 0$. Using this information we conclude that, there exists $C>0$ such that
\begin{equation*}
	\begin{aligned}
		J(t\phi) &= \frac{t^2}{2}\int\limits_{\RR^N} |\nabla \phi|^2 + a(x)\phi^2dx  -\int\limits_{\RR^N}F(x,t\phi)dx -\frac{1}{2\al\cdot{2^*_{\mu } }}\|G^{-1}(t\phi)\|_{0,\al}^{2\al\cdot2^{*}_{\mu }}\\
		&\leq \frac{Ct^2}{2}\| \phi\|^2 -\frac{\left( G^{-1}(t)\right)^{2\al\cdot{2^*_{\mu }}} }{2\al\cdot{2^*_{\mu } }}\|\phi\|_{0,\al}^{2\al\cdot2^{*}_{\mu }} \to -\infty \; \text{as} \; t \to \infty.
	\end{aligned}
\end{equation*}
Therefore, for sufficiently large $t$, we have $e = t\phi \in H^1(\RR^N)$, with $\|e\| > \rho$, such that $J(e) < 0$. Hence, the proof follows. \qed
\end{lem}
Let us consider the mini-max value
\begin{equation}\label{Deq3.3n}
	c := \inf\limits_{h \in \Gamma}\max\limits_{t \in [0,1]}J(h(t)),
\end{equation}
where 
\begin{equation*}
	\Gamma = \{h \in C([0,1], H^1(\RR^N)): h(0) = 0\;\text{and}\; h(1) = e\}.
\end{equation*} 
As a result of Lemma \ref{Dlem3.7}, there exists a $(PS)_c$ sequence $\{v_n\}$ in $H^1(\RR^N )$ at the level $c$, that is,
\begin{equation*}
	J(v_n) \to c \quad \text{and}\quad J'(v_n) \to 0 \;\text{in}\; H^{-1}(\RR^N ) \;\text{as}\; n \to +\infty.
\end{equation*}
Then, we have the following
\begin{lem}\label{Dlem3.8}
	Let  $\{v_n\}$ be a $(PS)_c$ sequence for the functional $J$, then $\{v_n\}$ is a bounded sequence in $H^1(\RR^N)$.
\end{lem}
\proof Let  $\{v_n\} \subset H^1(\RR^N)$ be a $(PS)_c$ sequence. Then, as $n \to \infty$, we have
\begin{equation*}
	J(v_n) \rightarrow c \;\; \text{and}\;\; |\langle J'(v_n), \phi\rangle|= o(1) \|\phi\|\;\; \text{for all} \;\phi \in C^\infty_0(\RR^N), \;\text{i.e.}
\end{equation*}
Employing Lemma \ref{Dlem3.1} (2), (4) and the fact that $2 < \tilde{\mu} < 2^* < 2\cdot2^*_\mu$, we deduce that
\begin{equation}\label{Deq3.5}
	\begin{aligned}
		\tilde{\mu}c& +o(1)
		 \geq\, \tilde{\mu}J(v_n)- \langle J'(v_n), v_n\rangle
		 \geq  \,\frac{\tilde{\mu}-2}{2}\int\limits_{\RR^N} |\nabla v_n|^2 +  a(x)|G^{-1}(v_n)|^2\,dx.
	\end{aligned}
\end{equation}
We divide the proof into the following cases:\\
\textbf{Case (i)} When $G^{-1}(v_n) > 1$,\\
Using Lemma \ref{Dlem3.1}(5), we have
\begin{equation*}
	H(x,t) \geq C(G(t))^{\tilde{\mu}} \geq C(G(t))^2\;\; \text{for all}\;\; x \in \RR^N , t > 1,
\end{equation*}
and from the fact that $J(v_n) = c + o(1)$ as $n \to \infty$, we infer that
\begin{equation*}
	\begin{aligned}
		\int\limits_{\{x:\;|G^{-1}(v_n)| > 1\}} a(x)|v_n|^2\,dx 
		\leq & C \int\limits_{\RR^N} H(x, G^{-1}(v_n))\,dx + \frac{C}{2\al\cdot2^*_{\mu}}\|G^{-1}(v^+_n)\|_{0,\al}^{2\al\cdot2^*_{\mu}}\\
		\leq & C\left[ \frac{1}{2}\int\limits_{\RR^{N}}|\nabla v_n|^2\,dx + \frac{1}{2}\int\limits_{\RR^{N}}a(x)|G^{-1} (v_n)|^2\,dx -c +o(1)\right].
	\end{aligned}
\end{equation*}
\textbf{Case (ii)} When $G^{-1}(v_n) \leq 1$,\\
As $v_n \leq G^{-1}(v_n)g(1)$, we get\\
\begin{equation*}
	\int\limits_{\{x:\;|G^{-1}(v_n)| \leq 1\}} \frac{a(x)|v_n|^2}{g^2(1)}\,dx 
	\leq \int\limits_{\RR^N} a(x)| G^{-1}(v_n)|^2\,dx.
\end{equation*}
Thus from both cases, we conclude that for some positive constant $C_1$
\begin{equation}\label{Deq3.6}
	\|v_n\|^2 = \int\limits_{\RR^{N}}|\nabla v_n|^2 + a(x)|v_n|^2\,dx
	\leq  C_1\left( \int\limits_{\RR^N} |\nabla v_n|^2 +  a(x)|G^{-1}(v_n)|^2\,dx\right) -c\cdot C + o(1). 
\end{equation}
By \eqref{Deq3.5} and \eqref{Deq3.6}, we deduce that
\begin{equation*}
	\|v_n\|^2 \leq \left[ \frac{2\tilde{\mu}C_1}{\tilde{\mu}-2}-C\right] c + o(1).
\end{equation*}
This implies that $\{v_n\}$ is bounded in $H^1(\RR^N)$.\qed

\section{Limiting equation at infinity}

In this section, by the Mountain-Pass Lemma and the second concentration-compactness principle, we prove the existence of a positive radial solution for $(P^\infty)$ which is the limit equation of \eqref{Deq1.6} at infinity.\\
 Proceeding as in Lemma \ref{Dlem3.7}, it is easy to verify the following: 
\begin{lem}\label{Dlem4.1} The functional $J^\infty(v)$ exhibits the Mountain-Pass geometry.
\end{lem}
Let $c^\infty$ be the mini-max value
\begin{equation*}
	c^\infty := \inf\limits_{h \in \Gamma}\max\limits_{t \in [0,1]}J^\infty(h(t)),
\end{equation*}
where 
\begin{equation*}
	\Gamma^\infty = \{h \in C([0,1], H^1(\RR^N)): h(0) = 0\;\text{and}\; J^\infty(h(1))< 0\}.
\end{equation*}
From Lemma \ref{Dlem4.1}, it follows that there exists a $(PS)_{c^\infty}$ sequence $\{v_n\}$ in $H^1(\RR^N )$ at the level $c^\infty$,
\begin{equation*}
	J^\infty(v_n) \to c^\infty \quad \text{and}\quad (J^\infty)'(v_n) \to 0 \;\text{in}\; H^{-1}(\RR^N ) \;\text{as}\; n \to +\infty.
\end{equation*}
By claiming as in Lemma \ref{Dlem3.8}, we obtain the boundedness of Palais-Smale sequences:
\begin{lem}\label{Dlem4.2}
	Let  $\{v_n\}$ be a $(PS)_{c^\infty}$ sequence for the functional $J^\infty$, then $\{v_n\}$ is a bounded sequence in $H^1(\RR^N)$.
\end{lem}
Next by Concentration-compactness principle, we prove the following compactness lemma.
\begin{lem}\label{Dprop4.1}
	Let $\{v_n\}$ be a $(PS)_{c^\infty}$ for $J^\infty$ defined on $H^1_r(\RR^N)$ with 
	\begin{equation*}
		c^\infty <	c_{*}^\infty :=\frac{1}{\al}\left[ \frac{1}{2}- \frac{1}{2\cdot2^*_{\mu}}\right]\left( \frac{\ba^2 S^H}{\al}\right)^{\frac{2^*_{\mu}}{2^*_{\mu}-1}}.
	\end{equation*}
	Then $\{v_n\}$ contains a convergent sub-sequence.
\proof Let $\{v_n\}$ be a $(PS)_{c^\infty}$ sequence for $J^\infty$, then by Lemma \ref{Dlem4.2}, we have $\{v_n\}$ is a bounded sequence thus $\{G^{-1}(v_n)\}$ is also bounded. Therefore there exists $v \in H^1(\RR^N)$ such that up to a sub-sequence $v_n \rightharpoonup v$ weakly in $H^1(\RR^N)$, $v_n \to v$ a.e. in $\RR^N$ and from the continuity of $G^{-1}$, we have the following
$$G^{-1}(v_n) \rightharpoonup G^{-1}(v) \;\text{ weakly in} \; H^1(\RR^N)\; \text{and}\; G^{-1}(v_n) \to G^{-1}(v) \;\text{a.e. in} \; \RR^N.$$
Furthermore, there exist bounded non-negative Radon measures $\omega$, $\xi$, and $\nu$ such that as $n \to \infty$
\begin{equation}\label{Deq4.2}
	|\nabla(G^{-1})^\al (v_n)|^2\rightharpoonup \omega, \;|G^{-1} (v_n)|^{\al\cdot2^*} \rightharpoonup \xi \;\;\text{and}\;\; \left(I_{N-\mu} *|G^{-1}(v_n^+)|^{\al\cdot2^*_{\mu}}\right)|G^{-1}(v_n^+)|^{\al\cdot2^*_{\mu}} \rightharpoonup \nu,
\end{equation} weakly in the sense of measures. Hence, by \cite[Lemma 2.5]{choq}
there exist an at most countable set $I$, a sequence of distinct points $\{x_i\}_{i\in I} \subset \RR^N$ and family of positive numbers $\{\nu_i\}_{i\in I}$, $\{\omega_i\}_{i\in I}$ and $\{\xi_i\}_{i\in I}$ such that
\begin{equation}\label{Deq4.3}
	\begin{aligned}
		\nu &= (I_{N-\mu} * |G^{-1}(v^+)|^{\al\cdot2^*_{\mu}})|G^{-1}(v^+)|^{\al\cdot2^*_{\mu}}+ \sum\limits_{i \in I}\nu_i\de_{x_i}, \;\; \sum\limits_{i \in I}\nu_i^{\frac{N}{2N-\mu}} < \infty;\\
		\xi &\geq |G^{-1}(v)|^{\al\cdot2^*} + \sum\limits_{i \in I}\xi_i\de_{x_i},\;\; \xi_i \geq C(N, \mu)^\frac{-N}{2N-\mu}\nu_{i}^\frac{N}{2N-\mu}\;\; \text{and}\\
		\omega &\geq |\nabla (G^{-1})^\al(v)|^2+ \sum\limits_{i \in I}\omega_i\de_{x_i},\;\; \omega_i \geq S^H\nu_i^{\frac{1}{2^*_{\mu}}},
	\end{aligned}
\end{equation}
where $\de_x$ is the Dirac-mass of mass 1 concentrated at $x \in \RR^N$.
Let $\e> 0$, we fix a smooth cut-off function $\phi_{\e,i}$ centred at $x_i$ such that 
\begin{equation*}
	0 \leq \phi_{\e,i} \leq 1,\; \phi_{\e,i} \equiv 1 \;\text{in}\; B(x_i, \e/2),\;  \phi_{\e,i} \equiv 0 \;\text{in}\; \RR^N\backslash B(x_i, \e) \;\text{and}\; |\nabla \phi_{\e,i}(x) | \leq \frac{4}{\e}.
\end{equation*}
Let $w_n = g(G^{-1}(v_n)) G^{-1}(v_n)$, then it is easy to verify using Lemma \ref{Dlem3.1} and assumption $(g_1)$(a) that $w_n$ is bounded in $H^1(\RR^N)$.
Thus, by dominated convergence theorem, we have
\begin{equation*}
	\begin{aligned}
		&\int\limits_{\RR^N}|G^{-1}(v_n)|^2\phi_{\e,i}(x)\;dx \to \int\limits_{\RR^N}|G^{-1}(v)|^2\phi_{\e,i}(x)\;dx,\\ &\int\limits_{\RR^N}\overline{h}(G^{-1}(v_n))G^{-1}(v_n)\phi_{\e,i}(x)\;dx  \to \int\limits_{\RR^N}\overline{h}(G^{-1}(v))G^{-1}(v)\phi_{\e,i}(x)\;dx, \;\text{as}\; n \to \infty
	\end{aligned}
\end{equation*}
and as $\e \to 0$ we have
\begin{equation}\label{Deq4.4}
	\begin{aligned}
		&\lim\limits_{\e \to 0}\lim\limits_{n \to \infty}\int\limits_{\RR^N}|G^{-1}(v_n)|^2\phi_{\e,i}(x)\;dx = 0,\\
		&\lim\limits_{\e \to 0}\lim\limits_{n \to \infty}\int\limits_{\RR^N}\overline{h}(G^{-1}(v_n))G^{-1}(v_n)\phi_{\e,i}(x)\;dx = 0.
	\end{aligned}
\end{equation}
In addition to this, employing H\"older's inequality, we get for some positive constant C
\begin{equation*}
	\int\limits_{\RR^N}\nabla v_n\cdot \nabla\phi_{\e,i}w_n \,dx \leq C \left(\int\limits_{\RR^N} |\nabla\phi_{\e,i}w_n|^2 \,dx \right)^{\frac{1}{2}}\leq \al C \left(\int\limits_{\RR^N} |\nabla\phi_{\e,i}v_n|^2 \,dx \right)^{\frac{1}{2}}.
\end{equation*}
Hence 
\begin{equation}\label{Deq4.5}
	\lim\limits_{\e \to 0}\lim\limits_{n \to \infty}\int\limits_{\RR^N}\nabla v_n\cdot \nabla\phi_{\e,i}w_n \,dx  = 0.
\end{equation}
Consider
\begin{equation}\label{Deq4.6}
\begin{aligned}
	& \left\langle (J^\infty)'(v_n), w_n\phi_{\e,i} \right\rangle 
	= \int\limits_{\RR^{N}}\nabla v_n\cdot(\nabla w_n \phi_{\e,i} +\nabla\phi_{\e,i}w_n)\,dx +\int\limits_{\RR^N}|G^{-1}(v_n)|^2\phi_{\e,i}(x)\;dx\\ 
	 \quad & - \int\limits_{\RR^N}\overline{h}(G^{-1}(v_n))G^{-1}(v_n)\phi_{\e,i}(x)\;dx  -\iint\limits_{\RR^{2N}}\frac{|G^{-1}(v_n^+(y))|^{\al\cdot2^*_{\mu}}|G^{-1}(v_n^+(x))|^{\al\cdot2^*_{\mu}}\phi_{\e,i}(x)}{|x-y|^\mu}\,dxdy.
	\end{aligned}
\end{equation}
Further using $(g_1)(b)$ and from Remark \ref{Drem1.1}, we deduce that
\begin{equation}\label{Deq4.7}
	\begin{aligned}
		\frac{g'(G^{-1}(v_n))}{g(G^{-1}(v_n))}G^{-1}(v_n) &\geq \ba^2 (\al-1)\frac{(G^{-1}(v_n))^{2\al-2}}{\left(g(G^{-1}(v_n))\right)^2};\\
		|\nabla v_n|^2 &\geq \ba^2 \frac{(G^{-1}(v_n))^{2\al-2}}{(g(G^{-1}(v_n)))^2}|\nabla v_n|^2.
	\end{aligned}
\end{equation}
Thus by \eqref{Deq4.7}, we get
\begin{equation}\label{Deq4.8}
	\begin{aligned}
		\nabla v_n \cdot \nabla w_n &= |\nabla v_n|^2 + \frac{g'(G^{-1}(v_n))}{g(G^{-1}(v_n))}G^{-1}(v_n)|\nabla v_n|^2\\
		&
		= \frac{\ba^2}{\al}|\nabla (G^{-1})^\al(v_n)|^2.
\end{aligned}\end{equation}
Using \eqref{Deq4.8}, the weak convergence of measure \eqref{Deq4.2} and \eqref{Deq4.3}, we deduce that
\begin{equation}\label{Deq4.9}
	\begin{aligned}
		&\lim\limits_{\e \to 0}\lim\limits_{n \to \infty}\int\limits_{\RR^{N}}\nabla v_n\cdot\nabla w_n \phi_{\e,i} \geq \frac{\ba^2}{\al} \lim\limits_{\e \to 0}\int\limits_{\RR^{N}}\phi_{\e,i}(x)d\omega \geq \frac{\ba^2}{\al} \omega_i.\\
		&\lim\limits_{\e \to 0}\lim\limits_{n \to \infty}\iint\limits_{\RR^{2N}}\frac{|G^{-1}(v_n^+(y))|^{\al\cdot2^*_{\mu}}|G^{-1}(v_n^+(x))|^{\al\cdot2^*_{\mu}}\phi_{\e,i}(x)}{|x-y|^\mu}\,dxdy = \lim\limits_{\e \to 0}\iint\limits_{\RR^{2N}}\phi_{\e,i}(x)d\nu = \nu_i.
	\end{aligned}
\end{equation}
Putting together \eqref{Deq4.4}, \eqref{Deq4.5} and \eqref{Deq4.9} in \eqref{Deq4.6}, we obtain
\begin{equation*}
	0 = \lim\limits_{\e \to 0}\lim\limits_{n \to \infty}\left\langle (J^\infty)'(v_n), w_n\phi_{\e,i} \right\rangle
	\geq \frac{\ba^2}{\al}\omega_i  - \nu_{i}.
\end{equation*}
It implies $\dfrac{\ba^2}{\al}\omega_i \leq \nu_{i}$. Combining this with the fact that $ S^H \nu_{i}^\frac{1}{2^*_{\mu}} \leq \omega_i$ we obtain
\begin{equation*}
	\omega_i \geq \left( \frac{\ba^2}{\al} (S^H)^{2^*_{\mu}} \right)^\frac{1}{2^*_{\mu}-1} \,\,\text{or}\;\; \omega_i = 0.
\end{equation*}
Let if possible, there exists $i_0 \in I$ such that $ \omega_{i_0} \geq \left( \dfrac{\ba^2}{\al} (S^H)^{2^*_{\mu}} \right)^\frac{1}{2^*_{\mu}-1} $. In particular, $\{v_n\}$ is a $(PS)_{c^\infty}$ sequence and employing $(g_1) (a)$ and Lemma \ref{Dlem3.1}(4), we deduce that
\begin{equation*}
	\begin{aligned}
		c^\infty &=\lim\limits_{\e \to 0}\lim\limits_{n \to \infty} \left( J^\infty(v_n)-\frac{1}{\tilde{\mu}\al}\left\langle (J^\infty)'(v_n), w_n \right\rangle\right)  \hspace{10cm}\\
		&\geq  \lim\limits_{\e \to 0}\lim\limits_{n \to \infty} \left(\left( \frac{1}{2}- \frac{1}{\tilde{\mu}}\right) \int\limits_{\RR^N}|\nabla v_n|^2\,dx+ \left( \frac{1}{\tilde{\mu}\al} -\frac{1}{2\al\cdot2^*_{\mu}}\right) \|G^{-1}(v_n^+)\|_{0,\al}^{2\al\cdot2^*_{\mu}}\right)\\
		&\geq \lim\limits_{\e \to 0}\lim\limits_{n \to \infty} \left\lbrace \left( \frac{1}{2}- \frac{1}{\tilde{\mu}}\right) \frac{\ba^2}{\al^2}\int\limits_{\RR^N}|\nabla (G^{-1})^\al(v_n)|^2\phi_{\e,i}\,dx \right.\\
		& \quad+ \left.\left( \frac{1}{\tilde{\mu}\al} -\frac{1}{2\al\cdot2^*_{\mu}}\right) \iint\limits_{\RR^{2N}}\frac{|G^{-1}(v_n^+(y))|^{\al\cdot2^*_{\mu}}|G^{-1}(v_n^+(x))|^{\al\cdot2^*_{\mu}}\phi_{\e,i}(x)}{|x-y|^\mu}\,dxdy\right\rbrace\\
		&\geq  \left( \frac{1}{2} - \frac{1}{\tilde{\mu}}\right)\frac{\ba^2}{\al^2}\omega_i + \left( \frac{1}{\tilde{\mu}\al} -\frac{1}{2\al\cdot2^*_{\mu}}\right)\nu_i\\
		&\geq  \left( \frac{1}{2}-\frac{1}{2\cdot2^*_{\mu}}\right)\frac{\ba^2}{\al^2}\omega_i
		= \frac{1}{\al}\left[ \frac{1}{2}-\frac{1}{2\cdot2^*_{\mu}}\right]\left(\frac{\ba^2}{\al}S^H\right)^{\frac{2^*_{\mu}}{2^*_{\mu}-1}}  > c^\infty.
	\end{aligned}
\end{equation*}
Therefore, $\omega_i = 0$ for all $i \in I$. Hence, we get
$$\|G^{-1}(v_n)\|_{0,\al}^{2\al\cdot2^*_{\mu}} \to \|G^{-1}(v)\|_{0,\al}^{2\al\cdot2^*_{\mu}}\; \text{as} \; n \to \infty.$$
Taking into account $(J^\infty)'(v_n) \to 0$ and using Br\'ezis-Lieb
lemma \cite[Theorem 1]{BrezLieb} and Fatou's Lemma, we have
\begin{equation}\label{Deq4.11}
	\begin{aligned}
		o(1) &=  \left\langle (J^\infty)'(v_n), w_n \right\rangle = \int\limits_{\RR^{N}}\nabla v_n\cdot\nabla w_n \,dx +\int\limits_{\RR^{N}}|v_n|^2\,dx - \int\limits_{\RR^{N}}k(v_n)\,dx - \|G^{-1}(v_n)\|_{0,\al}^{2\al\cdot2^*_{\mu}}\\
		&= \|v_n\|^2 + \int\limits_{\RR^{N}}\frac{g'(G^{-1}(v_n))}{g(G^{-1}(v_n))}G^{-1}(v_n)|\nabla v_n|^2\,dx - \int\limits_{\RR^{N}}k(v_n)\,dx - \|G^{-1}(v)\|_{0,\al}^{2\al\cdot2^*_{\mu}}\\
		&\geq \|v_n-v\|^2 + \|v\|^2 + \int\limits_{\RR^{N}}\frac{g'(G^{-1}(v))}{g(G^{-1}(v))}G^{-1}(v)|\nabla v|^2\,dx - \int\limits_{\RR^{N}}k(v_n)\,dx - \|G^{-1}(v)\|_{0,\al}^{2\al\cdot2^*_{\mu}},
	\end{aligned}
\end{equation}
where $k(t)= t^2 - |G^{-1}(t)|^2 + \overline{h}(G^{-1}(t))G^{-1}(t)$ and it has the following growth condition:
$$\lim\limits_{t \to 0} \dfrac{k(t)}{t^2} = 0;\;\text{and}\; \lim\limits_{t \to \infty} \dfrac{k(t)}{t^{2^*}} = 0.$$
Thus from \eqref{Deq4.11}, we get
\begin{equation*}
	o(1) \geq \|v_n-v\|^2 + \left\langle (J^\infty)'(v), w \right\rangle = \|v_n-v\|^2,
\end{equation*}
where $w = g(G^{-1}(v)) G^{-1}(v)$. Hence $v_n \to v $ in $H^1(\RR^N)$, finishing the proof.\qed
\end{lem}

\begin{lem}\label{Dlem4.3} Under the assumptions $(g_0)-(g_1)$, $(h_0)-(h_3)$ and $(a_0)$, we claim that
	$$\sup\limits_{t\geq 0} J^\infty(tu_\e)< c_{*}^\infty =\frac{1}{\al}\left[ \frac{1}{2}- \frac{1}{2\cdot2^*_{\mu}}\right]\left( \frac{\ba^2 S^H}{\al}\right)^{\frac{2^*_{\mu}}{2^*_{\mu}-1}}.$$
	\proof By Lemma \ref{Dlem4.1} we see that $J^\infty(0) = 0$ and $\lim\limits_{t \to \infty}J^{\infty}(tu_\e) = -\infty$. Thus for $\e>0$, there exists $t_\e>0$ such that
	$$\sup\limits_{t\geq 0}J^{\infty}(tu_\e) = J^{\infty}(t_\e u_\e).$$
	Also we claim that there exist positive constants $A_1$ and $A_2$ independent of $\e$ such that $0 < A_1 < t_\e < A_2 < \infty$. Next we denote $\ds\int\limits_{\RR^N}F^{*}(t_\e u_\e):= \|G^{-1}(t_\e u_\e)\|_{0,\al}^{2\al\cdot2^{*}_{\mu }} - \left( \frac{\al}{\ba}\right)^{2\cdot2^*_{\mu}}\|t_\e u_\e\|_0^{2\cdot2^*_{\mu}}$. Then we have
	\begin{equation}\label{Deq4.12}
			\sup\limits_{t\geq 0}J^{\infty}(tu_\e) = J^\infty(t_\e u_\e)
			 \leq \sup\limits_{t\geq 0}K(t) + \frac{t_\e^2}{2}\int\limits_{\RR^N} |u_\e|^2dx -\int\limits_{\RR^N}\overline{F}(t_\e u_\e)-\frac{1}{2\al\cdot{2^*_{\mu } }}\int\limits_{\RR^N}F^{*}(t_\e u_\e),
	\end{equation}
where \begin{equation*}
	K(t) = \ds \frac{t^2}{2}\int\limits_{\RR^N} |\nabla u_\e|^2 -\frac{t^{2\cdot2^*_{\mu}}}{2\al\cdot2^*_{\mu}}\left( \frac{\al}{\ba}\right)^{2\cdot2^*_{\mu}}\| u_\e\|_0^{2\cdot2^*_{\mu}}.
\end{equation*}
We observe that, $K(t)\to -\infty$ as $t\to\infty$ and $K(t)>0$ for small $t$.
So, there exists $t_{max}>0$ such that $\sup\limits_{t\geq 0}K(t)=K(t_{max})$. \begin{equation*}
	t_{max} = \left[ \frac{\ba^{2\cdot2^*_{\mu}}}{\al^{2\cdot2^*_{\mu}-1}} \frac{\|\nabla u_\e\|^2_{L^2}}{\| u_\e\|_0^{2\cdot2^*_{\mu}}} \right]^{\frac{1}{2\cdot2^*_{\mu}-2}}, \quad \;
	K(t_{max}) = \frac{1}{\al}\left( \frac{\ba^2}{\al}\right) ^{\frac{2^*_{\mu}}{2^*_{\mu}-1}}\left[ \frac{1}{2}- \frac{1}{2\cdot2^*_{\mu}}\right] \left[\frac{\|\nabla u_\e\|^2_{L^2}}{\| u_\e\|_0^{2}} \right]^{\frac{2^*_{\mu}}{2^*_{\mu}-1}}.
\end{equation*}
From Proposition \ref{Dprop2.2}, we have
\[\|\nabla u_\e\|^ {2}_{L^2}= C^\frac{N}{2\cdot2^*_{\mu}}(S^H)^\frac{N}{2}+ O(\e^{N-2})\; \text{and}
\;\; \|u_\e\|_0^{2}\geq {\left[C^\frac{N}{2}(S^H)^\frac{2N-\mu}{2}- O(\e^{N})\right] }^\frac{1}{2^*_{\mu}}.\]
This implies
for $\e$ small enough,
\begin{equation}\label{Deq4.13}
K(t_{max})\leq\frac{1}{\al}\left[ \frac{1}{2}- \frac{1}{2\cdot2^*_{\mu}}\right] \left(\frac{\ba^2}{\al}S^H\right)^{\frac{2^*_{\mu}}{2^*_{\mu}-1}} + O(\e^{N-2}) = c^\infty_{*}+ O(\e^{N-2}).
\end{equation}
 Further using Lemma \ref{Dlem3.1} and \eqref{Deq2.4} we see
 \begin{equation}\label{Deq4.14}
 	\begin{aligned}
 		\int\limits_{\RR^{N}}\overline{F}(t_\e u_\e)\,dx 
 		\geq \int\limits_{B_{2\rho_\e}}\overline{H}(G^{-1}(t_\e u_\e))\,dx
 		\geq t_\e^{\tilde{\mu}}\int\limits_{B_{2\rho_\e}}u_\e^{\tilde{\mu}}\,dx \geq A_1^{\tilde{\mu}}\int\limits_{B_{\rho_\e}}U_\e^{\tilde{\mu}}\,dx.
 		\end{aligned}
 \end{equation}
Note that, for some positive constant C, 
\begin{equation}\label{Deq4.15}
	\int\limits_{B_{\rho_\e}}U_\e^{\tilde{\mu}}\,dx = C\e^{N-\frac{(N-2)\tilde{\mu}}{2}}\int_{0}^{\frac{\rho_\e}{\e}}\frac{r^{N-1}}{(1+r^2)^{\frac{(N-2)\tilde{\mu}}{2}}}\,dr \geq C\e^{N-\frac{(N-2)\tilde{\mu}}{2}}.
\end{equation}
Putting together \eqref{Deq4.13}, \eqref{Deq4.14} and \eqref{Deq4.15} in \eqref{Deq4.12} 
		\begin{equation}\label{Deq4.16}
		\sup\limits_{t\geq 0}J^{\infty}(tu_\e)
			 \leq c^\infty_{*}+ O(\e^{N-2}) + \frac{t_\e^2}{2}\int\limits_{B_{2\rho_\e}} |u_\e|^2dx -A_1^{\tilde{\mu}}C\e^{N-\frac{(N-2)\tilde{\mu}}{2}}-\frac{1}{2\al\cdot{2^*_{\mu } }}\int\limits_{B_{2\rho_\e}}F^{*}(t_\e u_\e).
	\end{equation}
Now, we are going to estimate $F^{*}(t_\e u_\e)$.	\\
\textbf{Case I:} Evaluating on $B_{\rho_\e} \times B_{\rho_\e}$\\
Using Lemma \ref{Dlem3.1}(6), the symmetry of variables and Hardy-Littlewood-Sobolev inequality, we get for some positive constant $C$
\begin{equation}\label{Deq4.17}
	\begin{aligned}
		&\iint\limits_{B_{\rho_\e} \times B_{\rho_\e}}\frac{\left( G^{-1}(t_\e u_\e(y))\right) ^{\al\cdot2^{*}_{\mu }}\left( G^{-1}(t_\e u_\e(x))\right) ^{\al\cdot2^{*}_{\mu }}- \left( \frac{\al}{\ba}\right)^{2\cdot2^*_{\mu}}\left( t_\e u_\e(y)\right) ^{2^*_{\mu}}\left( t_\e u_\e(x)\right) ^{2^*_{\mu}}}{|x-y|^{\mu}}\,dxdy\\
		&= \iint\limits_{B_{\rho_\e} \times B_{\rho_\e}}\frac{\left( G^{-1}(t_\e u_\e(y))\right) ^{\al\cdot2^{*}_{\mu }}\left( t_\e u_\e(x)\right) ^{2^*_{\mu}}}{|x-y|^{\mu}}\left[\frac{\left( G^{-1}(t_\e u_\e(x))\right) ^{\al\cdot2^{*}_{\mu }}}{\left( t_\e u_\e(x)\right) ^{2^*_{\mu}}}- \left( \frac{\al}{\ba}\right)^{2^*_{\mu}} \right] \,dxdy\\
		 &\quad + \left( \frac{\al}{\ba}\right)^{2^*_{\mu}}\iint\limits_{B_{\rho_\e} \times B_{\rho_\e}}\left[\frac{\left( G^{-1}(t_\e u_\e(y))\right) ^{\al\cdot2^{*}_{\mu }}}{\left( t_\e u_\e(y)\right) ^{2^*_{\mu}}}- \left( \frac{\al}{\ba}\right)^{2^*_{\mu}} \right]\frac{\left( t_\e u_\e(y)\right) ^{2^*_{\mu}}\left( t_\e u_\e(x)\right) ^{2^*_{\mu}}}{|x-y|^{\mu}}\,dxdy\\
		 &\geq -C\iint\limits_{B_{\rho_\e} \times B_{\rho_\e}}\frac{\left( G^{-1}(t_\e u_\e(y))\right) ^{\al\cdot2^{*}_{\mu }}\left( t_\e u_\e(x)\right) ^{2^*_{\mu}-\de}}{|x-y|^{\mu}}+ \left( \frac{\al}{\ba}\right)^{2^*_{\mu}}\frac{\left( t_\e u_\e(y)\right) ^{2^*_{\mu}-\de}\left( t_\e u_\e(x)\right) ^{2^*_{\mu}}}{|x-y|^{\mu}}\,dxdy\\
		 & \geq -2C\left( \frac{\al}{\ba}\right)^{2^*_{\mu}}t_\e^{2\cdot2^*_{\mu}-\de}\iint\limits_{B_{\rho_\e} \times B_{\rho_\e}}\frac{\left(u_\e(y)\right) ^{2^*_{\mu}-\de}\left( u_\e(x)\right) ^{2^*_{\mu}}}{|x-y|^{\mu}}\,dxdy.\\
	&\geq -C \|u_\e\|_{L^{2^*}(B_{\rho_\e})}^{2^*_{\mu}} \left( \int\limits_{B_{\rho_\e}}(u_\e(y))^{\frac{(2^*_{\mu}-\de)2^*}{2^*_{\mu}}}\,dy\right)^{\frac{2^*_{\mu}}{2^*}}.
	\end{aligned}
	\end{equation} 		
		As $2^*_{\mu} > 2$, we have $t := \dfrac{(2^*_{\mu}-\de)2^*}{2^*_{\mu}} > \dfrac{N}{N-2}$, therefore by Proposition \ref{Dprop2.2}, there exists $C >0$ such that
		\begin{equation}\label{Deq4.18}
			\left(\int\limits_{B_{\rho_\e}}u_\e(y)^t\,dy\right)^{\frac{2^*_{\mu}}{2^*}} = C\e^{\frac{(N-2)\de}{2}} \left(1+ O\left(\e^{\frac{(2^*_{\mu}-\de)N}{2^*_{\mu}}}\right) \right).
		\end{equation} 
	Thus substituting \eqref{Deq4.18} in \eqref{Deq4.17} and using $\|u_\e\|_{L^{2^*}}$ estimate from Proposition \ref{Dprop2.2}, we get
		\begin{equation}\label{Deq4.19}
		\begin{aligned}
			-\int F^{*}(t_\e u_\e) &\leq C \left[S^{\frac{N}{2}} +O(\e^N) \right] ^{\frac{2^*_{\mu}}{2^*}} \e^{\frac{(N-2)\de}{2}} \left(1+ O(\e^{\frac{(2^*_{\mu}-\de)N}{2^*_{\mu}}}) \right)\\
			&= C\e^{\frac{(N-2)\de}{2}} \left(1+ O\left(\e^{\frac{(2^*_{\mu}-\de)N}{2^*_{\mu}}}\right) \right).
		\end{aligned} 
	\end{equation}
\textbf{Case II:} Evaluating on $D \times D$, where $D:= B_{2\rho_\e}\backslash B_{\rho_\e}$\\
Using Lemma \ref{Dlem3.1}(6), symmetry of variables and Hardy-Littlewood-Sobolev inequality, we get for some positive constant C
\begin{equation}\label{Deq4.20}
	\begin{aligned}
		&\iint\limits_{D \times D}\frac{\left( G^{-1}(t_\e u_\e(y))\right) ^{\al\cdot2^{*}_{\mu }}\left( G^{-1}(t_\e u_\e(x))\right) ^{\al\cdot2^{*}_{\mu }}- \left( \frac{\al}{\ba}\right)^{2\cdot2^*_{\mu}}\left( t_\e u_\e(y)\right) ^{2^*_{\mu}}\left( t_\e u_\e(x)\right) ^{2^*_{\mu}}}{|x-y|^{\mu}}\,dxdy\\
		&= \iint\limits_{D \times D}\frac{\left( G^{-1}(t_\e u_\e(y))\right) ^{\al\cdot2^{*}_{\mu }}\left( t_\e u_\e(x)\right) ^{2^*_{\mu}}}{|x-y|^{\mu}}\left[\frac{\left( G^{-1}(t_\e u_\e(x))\right) ^{\al\cdot2^{*}_{\mu }}}{\left( t_\e u_\e(x)\right) ^{2^*_{\mu}}}- \left( \frac{\al}{\ba}\right)^{2^*_{\mu}} \right] \,dxdy\\
		&\quad + \left( \frac{\al}{\ba}\right)^{2^*_{\mu}}\iint\limits_{D \times D}\left[\frac{\left( G^{-1}(t_\e u_\e(y))\right) ^{\al\cdot2^{*}_{\mu }}}{\left( t_\e u_\e(y)\right) ^{2^*_{\mu}}}- \left( \frac{\al}{\ba}\right)^{2^*_{\mu}} \right]\frac{\left( t_\e u_\e(y)\right) ^{2^*_{\mu}}\left( t_\e u_\e(x)\right) ^{2^*_{\mu}}}{|x-y|^{\mu}}\,dxdy\\
		&\geq -2\left( \frac{\al}{\ba}\right)^{2\cdot2^*_{\mu}}\iint\limits_{D\times D}\frac{\left( t_\e u_\e(y)\right) ^{2^*_{\mu}}\left( t_\e u_\e(x)\right) ^{2^*_{\mu}}}{|x-y|^{\mu}}\,dxdy\\
		& \geq -C\|u_\e\|_{L^{2^*}(D)}^{2\cdot2^*_{\mu}}.
	\end{aligned}
\end{equation}
We observe that, for $r \in (\rho_\e, 2\rho_e)$ there exist positive constants $C_1$ and $C_2$ such that
\begin{equation*}
	C_1\e^{\frac{N-2}{2}-(N-2)\tau} \leq U_\e(r) \leq C_2\e^{\frac{N-2}{2}-(N-2)\tau}.
\end{equation*}
Therefore using the above estimate in \eqref{Deq4.20}, we get
\begin{equation}\label{Deq4.21}
	-\int F^{*}(t_\e u_\e) \leq C \e^{2\cdot2^*_{\mu}\left( \frac{N-2}{2}-(N-2)\tau\right) } \left[ \int_{\rho_\e}^{2\rho_\e}r^{N-1}\right]^ {\frac{2\cdot2^*_{\mu}}{2^*}}
	= C\e^{(2N-\mu)(1-\tau)}.
\end{equation} 		
\textbf{Case III:} Evaluating on $B_{\rho_\e} \times D$, where $D:= B_{2\rho_\e}\backslash B_{\rho_\e}$\\
Similar to the cases above, employing Lemma \ref{Dlem3.1}(6) and symmetry of variables, we get for some positive constant C
\begin{equation}\label{Deq4.22}
	\begin{aligned}
		&\iint\limits_{B_{\rho_\e} \times D}\frac{\left( G^{-1}(t_\e u_\e(y))\right) ^{\al\cdot2^{*}_{\mu }}\left( G^{-1}(t_\e u_\e(x))\right) ^{\al\cdot2^{*}_{\mu }}- \left( \frac{\al}{\ba}\right)^{2\cdot2^*_{\mu}}\left( t_\e u_\e(y)\right) ^{2^*_{\mu}}\left( t_\e u_\e(x)\right) ^{2^*_{\mu}}}{|x-y|^{\mu}}\,dxdy\\
		&= \iint\limits_{B_{\rho_\e} \times D}\frac{\left( G^{-1}(t_\e u_\e(y))\right) ^{\al\cdot2^{*}_{\mu }}\left( t_\e u_\e(x)\right) ^{2^*_{\mu}}}{|x-y|^{\mu}}\left[\frac{\left( G^{-1}(t_\e u_\e(x))\right) ^{\al\cdot2^{*}_{\mu }}}{\left( t_\e u_\e(x)\right) ^{2^*_{\mu}}}- \left( \frac{\al}{\ba}\right)^{2^*_{\mu}} \right] \,dxdy\\
		&\quad + \left( \frac{\al}{\ba}\right)^{2^*_{\mu}}\iint\limits_{B_{\rho_\e} \times D}\left[\frac{\left( G^{-1}(t_\e u_\e(y))\right) ^{\al\cdot2^{*}_{\mu }}}{\left( t_\e u_\e(y)\right) ^{2^*_{\mu}}}- \left( \frac{\al}{\ba}\right)^{2^*_{\mu}} \right]\frac{\left( t_\e u_\e(y)\right) ^{2^*_{\mu}}\left( t_\e u_\e(x)\right) ^{2^*_{\mu}}}{|x-y|^{\mu}}\,dxdy\\
		&\geq -C\iint\limits_{B_{\rho_\e} \times D}\frac{\left( G^{-1}(t_\e u_\e(y))\right) ^{\al\cdot2^{*}_{\mu }}\left( t_\e u_\e(x)\right) ^{2^*_{\mu}-\de}}{|x-y|^{\mu}} + \left( \frac{\al}{\ba}\right)^{2\cdot2^*_{\mu}}\frac{t_\e^{2\cdot2^*_{\mu}}\left(  u_\e(y)\right) ^{2^*_{\mu}}\left( u_\e(x)\right) ^{2^*_{\mu}}}{|x-y|^{\mu}}\,dxdy\\
		& \geq -C\left( \frac{\al}{\ba}\right)^{2^*_{\mu}}t_\e^{2\cdot2^*_{\mu}-\de}\iint\limits_{B_{\rho_\e} \times D}\frac{\left(u_\e(y)\right) ^{2^*_{\mu}}\left( u_\e(x)\right) ^{2^*_{\mu}-\de}}{|x-y|^{\mu}}\,dxdy + \left( \frac{\al}{\ba}\right)^{2\cdot2^*_{\mu}}t_\e^{2\cdot2^*_{\mu}}\|u_\e\|_{0}^{2\cdot2^*_{\mu}}.
	\end{aligned}
\end{equation}
Further using Hardy-Littlewood-Sobolev inequality in \eqref{Deq4.22} and using the estimates from Case I and Case II, we imply that
\begin{equation}\label{Deq4.23}
\begin{aligned}
		-\int F^{*}(t_\e u_\e) &\leq C \|u_\e\|_{L^{2^*}(D)}^{2^*_{\mu}} \left( \int\limits_{B_{\rho_\e}}u_\e^{\frac{(2^*_{\mu}-\de)2^*}{2^*_{\mu}}}\right)^{\frac{2^*_{\mu}}{2^*}} + C\|u_\e\|_{L^{2^*}(B_{\rho_\e})}^{2^*_{\mu}}\|u_\e\|_{L^{2^*}(D)}^{2^*_{\mu}}\\
		&\leq C\e^{\frac{(2N-\mu)}{2}(1-\tau)+{\frac{(N-2)\de}{2}}} \left(1+ O\left(\e^{\frac{(2^*_{\mu}-\de)N}{2^*_{\mu}}}\right) \right) + C\e^{\frac{(2N-\mu)}{2}(1-\tau)}\left( 1 + O(\e^N) \right).
	\end{aligned}
\end{equation}
\textbf{Case IV:} Evaluating on $D \times B_{\rho_\e}$, where $D:= B_{2\rho_\e}\backslash B_{\rho_\e}$\\
It follows the same argument as in Case III.\\
Further assembling the estimates obtained in \eqref{Deq4.19}, \eqref{Deq4.21}, \eqref{Deq4.23} and substituting them in \eqref{Deq4.16}, we get for some positive constant C
\begin{equation}\label{Deq4.24}
	\begin{aligned}
		&\sup\limits_{t\geq 0}J^{\infty}(tu_\e)
		\leq c^\infty_{*}+ O(\e^{N-2}) -C\e^{N-\frac{(N-2)\tilde{\mu}}{2}} + C\e^{\frac{(N-2)\de}{2}}+  C\e^{(2N-\mu)(1-\tau)}
		+ C\e^{\frac{(2N-\mu)}{2}(1-\tau)}\\  
		&\quad+ C\e^{\frac{(2N-\mu)}{2}(1-\tau)+{\frac{(N-2)\de}{2}}} + \begin{cases}
			C\e^2+ O(\e^{N-2}) & N \geq 5,\\
			C\e^2|\ln(\e)|+ O(\e^2) & N= 4,\\
			O(\e) & N =3.\\
		\end{cases}
	\end{aligned}
\end{equation}
Thus by direct calculation, we conclude that if
\begin{enumerate}
\item [(1)] $N \geq \max\left\lbrace 2 + \frac{4\al}{\al(\tilde{\mu}-1)-\ga^+}, 6\right\rbrace;\, \tilde{\mu} >2$, 
\item [(2)] $N = 5;\, \tilde{\mu} > \frac{7}{3} +\frac{\ga^+}{\al},$
\item [(3)] $N = 4;\, \tilde{\mu} > 3 +\frac{\ga^+}{\al},$
\item [(4)] $N = 3;\, \tilde{\mu} > 5 +\frac{\ga^+}{\al},$
\end{enumerate}
then, we can choose $\tau \in ( \frac{1}{2} , 1)$, such that
\begin{equation*}
	\min\left\lbrace N-2, \frac{(N-2)\de}{2}, (2N-\mu)(1-\tau), \frac{(2N-\mu)}{2}(1-\tau)\right\rbrace >  N-\frac{(N-2)\tilde{\mu}}{2}.
\end{equation*}
Hence by \eqref{Deq4.24} and the above conditions on $N$ and $\tilde{\mu}$, we get that for $\e > 0$, sufficiently small
$$\sup\limits_{t\geq 0}J^{\infty}(tu_\e) < c^\infty_{*}.\qed$$ 
\end{lem}

\begin{thm}\label{Dthm4.1}
	Let  $(g_0)-(g_1)$, $(h_0)-(h_3)$ and $(a_0)$ hold. Then for $(P^\infty)$, there exists at least one positive radial solution provided $N$ and $\tilde{\mu}$ satisfies
	\begin{enumerate}
		\item[(1)] $N \geq \max\left\lbrace 2 + \frac{4\al}{\al(\tilde{\mu}-1)-\ga^+}, 6\right\rbrace$ and $ \tilde{\mu} > 2$,
		\item [(2)] $N = 5;\, \tilde{\mu} > \frac{7}{3} +\frac{\ga^+}{\al},$
		\item [(3)] $N = 4;\, \tilde{\mu} > 3 +\frac{\ga^+}{\al},$
		\item [(4)] $N = 3;\, \tilde{\mu} > 5 +\frac{\ga^+}{\al},$
	\end{enumerate} 
where $\ga^+ = \max \{\ga, 0\}$.
	\proof From Lemma \ref{Dlem4.1}, there exists a $(PS)_{c^\infty}$ sequence $\{v_n\}$ in $H^1(\RR^N )$ at the level $c^\infty$, that is,
	$J^\infty(v_n) \to c^\infty$ and $(J^\infty)'(v_n) \to 0$ in $H^{-1}(\RR^N )$ as $n \to +\infty$.
	Consequently, by Lemma \ref{Dprop4.1} and Lemma \ref{Dlem4.3} we get
	$0 < \mathfrak{a}  < c^\infty \leq \max\limits_{t \in [0,1]}J^\infty(tu_\e) < c^{\infty}_{*}$.
	Thus by the Mountain-pass theorem, there exists $v \in H^1(\RR^N)$ such that $v_n \to v$ in $H^1(\RR^N)$  i.e. $v$ is the non-trivial critical
	point of $J^\infty$. Moreover, by standard regularity argument, it follows that $v \in C^2(\RR^N)$ and by Lemma \ref{Dlem3.1} and the assumption on $h$, we have
	\begin{equation*}
		\begin{aligned}
			-\De v + v  &= \left( v-\frac{G^{-1}(v)}{g(G^{-1}(v))}\right)  +\frac{\overline{h}(G^{-1}(v))}{g(G^{-1}(v))} + \left( I_{\vartheta}*|G^{-1}(v^+)|^{\al\cdot2^*_\mu}\right) \frac{|G^{-1}(v^+)|^{\al\cdot2^*_\mu-2}G^{-1}(v^+)}{g(G^{-1}(v^+))}\\
			& \geq 0.
		\end{aligned}
	\end{equation*}
Thus from strong maximum principle, it follows that $v$ is a positive solution of problem $(P^\infty)$.\qed
\end{thm}

We define
$J^\infty = \inf\limits_{v \in \mathcal{M}}J^\infty(v)$,
where
\begin{equation*}
	\mathcal{M} = \{v \in H^1(\RR^N ) \backslash \{0\}; \left\langle (J^\infty)'(v), v\right\rangle = 0\}.
\end{equation*}
By Theorem \ref{Dthm4.1}, we conclude that $\mathcal{M} \neq \phi$.

\begin{lem}\label{Dlem4.4}
	Let  $(g_0)-(g_1)$, $(h_0)-(h_3)$ and $(a_0)$ hold, then $J^\infty = c^\infty$.
\proof Employing Theorem \ref{Dthm4.1}, we see that there exists $v \in \mathcal{M}$, such that $c^\infty = J^\infty(v) \geq \inf\limits_{u \in \mathcal{M}}J^\infty(u) = J^\infty$. \\
On the other hand, by the mountain-pass geometry, we know for all $u \in \mathcal{M}$, there exists $t_{*} > 0$, such that
\begin{equation*}
	J^\infty(t_{*}u) = \sup\limits_{t>0}J^\infty(tu)\; \;\text{and}\;\;\frac{d}{dt}{J^\infty(tu)\left|_{t=t_{*}}\right.}= 0. 
\end{equation*}
Denote
\begin{equation*}
\begin{aligned}
		k = \int\limits_{\RR^{N}} |\nabla u|^2\,dx; \quad Q_1(s) = \frac{G^{-1}(s)}{g(G^{-1}(s))}; \quad  
	Q_2(s) = \frac{\overline{h}(G^{-1}(s))}{g(G^{-1}(s))}.
	\end{aligned}	
\end{equation*}
By Lemma \ref{Dlem3.1}(2) and (4), we conclude that
\begin{equation}\label{Deq4.25}
	sQ'_1(s) - Q_1(s) \leq 0; \quad s\frac{d Q_2(s)}{ds}- Q_2(s) \geq 0 \quad \text{for all} \;s \geq 0.
\end{equation}
 Consider
\begin{equation*}
	\begin{aligned}
		\frac{dJ^\infty(tu)}{dt} &= tk + \int\limits_{\RR^N}Q_1(tu)u\,dx - \int\limits_{\RR^N}Q_2(tu)u\,dx - \iint\limits_{\RR^{2N}}\mathcal{A}(tu)u(x)\,dxdy\\
		&= t\nu(t)
	\end{aligned}
	\end{equation*}
 where $$\mathcal{A}(tu) = \ds \frac{|G^{-1}(tu(y))|^{\al\cdot2^*_\mu}|G^{-1}(tu(x))|^{\al\cdot2^*_\mu-2}G^{-1}(tu(x))}{|x-y|^{\mu}g(G^{-1}(tu(x)))}$$  and $$\nu(t) = \ds k + \int\limits_{\RR^N}\frac{Q_1(tu)u}{t}\,dx - \int\limits_{\RR^N}\frac{Q_2(tu)u}{t}\,dx - \iint\limits_{\RR^{2N}}\frac{\mathcal{A}(tu)}{t}u(x)\,dxdy.$$ 
 By direct calculation and using the assumption $(g_1) (a)$ 
 and Lemma \ref{Dlem3.1}(4)
 , we deduce that
 \begin{equation}\label{Deq4.26}
 	t\frac{d \mathcal{A}(tu)}{dt}- \mathcal{A}(tu) \geq 0.
 \end{equation}
 Employing \eqref{Deq4.25} and \eqref{Deq4.26}, we get for $t >0$
\begin{equation*}
	\begin{aligned}
		\nu'(t) &= \int\limits_{\RR^N}\frac{tu\frac{dQ_1(tu)}{dt}- Q_1(tu)}{t^2}u\,dx - \int\limits_{\RR^N}\frac{tu\frac{dQ_2(tu)}{dt}- Q_2(tu)}{t^2}u\,dx\\ 
		& \qquad- \iint\limits_{\RR^{2N}}\frac{t\frac{d\mathcal{A}(tu)}{dt}- \mathcal{A}(tu)}{t^2}u(x)\,dxdy\\
		& < 0.
	\end{aligned}
\end{equation*}
Hence $\nu(t)$ has at most one zero point in $(0,+\infty)$. We have assumed $\frac{dJ^\infty(tu)}{dt}\left|_{
	t=t_{*}} = 0\right.$ and since $u \in \mathcal{M},$ $\frac{dJ^\infty(tu)}{dt}\left|_{t=1} = 0\right.$, we can conclude that $t_{*} = 1$. Thus
$J^\infty= \inf\limits_{u \in \mathcal{M}}J^{\infty}(u) = \inf\limits_{u \in \mathcal{M}}\sup\limits_{t>0}J^\infty(tu) \geq  c^\infty.$ Thus we get the desired result. \qed
\end{lem}
From Lemma \ref{Dlem4.4}, we conclude that $J^\infty = c^\infty$ can be attained by a function $w\in \mathcal{M}$ and $w$ is a positive ground
state solution of
\begin{equation*}  
		-\De w + G^{-1}(w) =  \overline{h}(G^{-1}(w)) + \left(I_{\vartheta}*|G^{-1}(w)|^{\al\cdot2^*_\mu}\right)|G^{-1}(w)|^{\al\cdot2^*_\mu-2}G^{-1}(w).
\end{equation*}
In \cite{gidas.ni.nirenberg}, Gidas, Ni and Nirenberg proved that there exist $a_1$, $a_2 > 0$ such that,
\begin{equation}\label{Deq4.27}
	a_1(|x| + 1)^{\frac{-(N-1)}{2}}e^{-|x|}\leq  w(x) \leq a_2(|x| + 1)^{\frac{-(N-1)}{2}}e^{-|x|} \;\;\text{for all}\;\; x \in \RR^N.
\end{equation}
Using estimate \eqref{Deq4.27}, we can obtain the following lemma
\begin{lem}\label{Dlem4.5}\cite[Lemma 3.6]{chen}
	Let $\nu$ be a unit vector of $\RR^N$ and $w$ be that in \eqref{Deq4.27}. There exist some constants
	$C_1 > 0$, $C_2 > 0$, and $C_3 > 0$ independent of $R \geq 1$ such that
	
		\begin{enumerate}
			\item [(1)] $\ds\int\limits_{\{x\in \RR^N : |x|\leq1\}}(w(x - R\nu))^2\,dx \geq C_1R^{-(N-1)}e^{-2R}; $
			\item [(2)] $\ds\int\limits_{ \RR^N }e^{-\mu|x|}(w(x - R\nu))^2\,dx \leq C_2R^{-(N-1)}e^{-2R}; $
			\item [(3)] $\ds\int\limits_{ \RR^N }e^{-\mu|x|}(w(x - R\nu))^{p+1}\,dx \leq C_3e^{-\min\{\mu, p+1\}R}.$
		\end{enumerate}
\end{lem}

\section{Proof of Theorem \ref{Dthm1.1}}
This section is devoted to the proof of the Theorem \ref{Dthm1.1}. To obtain the existence of a positive solution, we deploy
the concentration-compactness principle and the classical mountain-pass theorem.
Moreover, from Lemma \ref{Dlem4.4} and Lemma \ref{Dlem4.3}, we can deduce that
\begin{equation}\label{Deq5.1n}
	J^\infty = c^\infty \leq \sup\limits_{t \geq 0}J^\infty(tu_\e) < c^\infty_{*}.
\end{equation}
Thus, it suffices to prove that $J$ satisfies Palais-Smale condition for $c \in \left(0,J^\infty\right)$.

\begin{lem}\label{Dlem5.1}
Assume  $(g_0)-(g_1)$, $(h_0)-(h_3)$ and $(a_0)$ hold. Let $\{v_n\}$ be a $(PS)_{c}$ for $J$ with $\ds c \in \left(0,J^\infty\right)$. Then $\{v_n\}$ contains a convergent sub-sequence.
\proof  Let $\{v_n\}$ be a $(PS)_c$ sequence for $J$ then by Lemma \ref{Dlem3.8}, we have that $\{v_n\}$ is bounded in $H^1(\RR^N )$.
Therefore there exists $v \in H^1(\RR^N)$ such that up to a sub-sequence 
\begin{equation}\label{Deq5.1}
	v_n \rightharpoonup v \;\text{weakly in}\; H^1(\RR^N),\;v_n \to v\;\text{in}\; L^s_{loc}(\RR^N)\;\text{for}\; s\in [1, 2^*),\; v_n \to v\;\text{a.e. in}\;\RR^N.
\end{equation}
Here $v$ is a weak solution of \eqref{Deq1.6} with $J(v) \geq 0$. Denote $w_n = v_n - v$. We will prove that $\|w_n\| \to 0$ as $n \to \infty$. Assume by contradiction that $\|w_n\| \to l > 0$ as $n \to \infty$.
Now we divide our proof into the following two steps:\\
\textbf{Step 1:} When $\|w_n\|_{L^2} \to 0$ as $n \to \infty$.\\
Thus by interpolation inequality, we claim that $\|w_n\|_{L^r} \to 0$ as $n \to \infty$ for $r \in [2, 2^*)$. In fact, since $2 <r< 2^*$, there exists $\te_r \in (0, 1)$ such that
\begin{equation*}
	\frac{1}{r} = \frac{\te_r}{2} + \frac{1-\te_r}{2^*}.
\end{equation*}
Since $\|w_n\|$ is bounded in $L^{2^*}$ norm, we get
\begin{equation*}
	\|w_n\|_{L^r} \leq \|w_n\|_{L^2}^{\te_r}\|w_n\|_{L^{2^*}}^{1-\te_r} \leq C\|w_n\|_{L^2}^{\te_r} \to 0 \;\text{as} \; n \to \infty.
\end{equation*}
Employing this claim and Lemma \ref{Dlem3.2}, we conclude that for any $\e > 0$, there exists $C_\e > 0$ such that
\begin{equation*}
	\left|\int\limits_{\RR^N}F(x,w_n)\,dx\right| \leq \e\left( \int\limits_{\RR^{N}} |w_n|^2 + |w_n|^{2^*}\,dx\right) + C_\e \int\limits_{\RR^{N}} |w_n|^r\,dx,
\end{equation*}
which implies that
\begin{equation}\label{Deq5.2}
	\int\limits_{\RR^N}F(x,w_n)\,dx =o(1), \;\text{as}\; n \to \infty.
\end{equation}
Likewise
\begin{equation}\label{Deq5.3}
	\int\limits_{\RR^N}f(x,w_n)\phi\,dx = o(1)\;\text{as}\; n \to \infty,\; \text{for any} \; \phi \in C^\infty_0(\RR^N),
\end{equation}
also we have
\begin{equation}\label{Deq5.4}
\left| \int\limits_{\RR^N}a(x)w_n^2\,dx\right|  \leq C \int\limits_{\RR^N} |w_n|^2\,dx  = o(1)\;\text{as}\; n \to \infty.
\end{equation}
Employing Brezis-Lieb \cite{BrezLieb, M.yang}, \cite[Lemma 2.2]{zhu_cao}, \eqref{Deq5.2} and \eqref{Deq5.4}, we get as $n \to \infty$
\begin{equation}\label{Deq5.5}
	\begin{aligned}
		c +o(1) &= J(v_n)\\
		& = \frac{1}{2}\int\limits_{\RR^N} |\nabla w_n|^2 dx -\frac{1}{2\al\cdot{2^*_{\mu } }}\|G^{-1}(w_n^+)\|_{0,\al}^{2\al\cdot2^{*}_{\mu }} + J(v).\\
	\end{aligned}
\end{equation}
By \eqref{Deq5.1}, \eqref{Deq5.3} and using the fact that $J'(v) = 0$, it follows that
\begin{equation}\label{Deq5.6}
	\begin{aligned}
		o(1) &= \left\langle J'(v_n) - J'(v), v_n -v\right\rangle \\
		& = \int\limits_{\RR^N} |\nabla w_n|^2 + a(x)w_n^2 dx -\iint\limits_{\RR^{2N}}\frac{|G^{-1}(v_n^+(y))|^{\al\cdot2^{*}_{\mu }}|G^{-1}(v_n^+(x))|^{\al\cdot2^{*}_{\mu}-2}G^{-1}(v_n^+(x))w_n}{|x-y|^ \mu \,g(G^{-1}(v_n^+(x)))}dxdy\\
		&\quad+ \iint\limits_{\RR^{2N}}\frac{|G^{-1}(v^+(y))|^{\al\cdot2^{*}_{\mu }}|G^{-1}(v^+(x))|^{\al\cdot2^{*}_{\mu}-2}G^{-1}(v^+(x))w_n}{|x-y|^ \mu \,g(G^{-1}(v^+(x)))}dxdy.
	\end{aligned}
\end{equation}
By reasoning as in Brezis-Lieb, we get
\begin{equation}\label{Deq5.7}
\begin{aligned}
		&\iint\limits_{\RR^{2N}}\frac{|G^{-1}(v_n^+(y))|^{\al\cdot2^{*}_{\mu }}|G^{-1}(v_n^+(x))|^{\al\cdot2^{*}_{\mu}-2}G^{-1}(v_n^+(x))w_n}{|x-y|^ \mu \,g(G^{-1}(v_n^+(x)))}dxdy \\ &=\iint\limits_{\RR^{2N}}\frac{|G^{-1}(v^+(y))|^{\al\cdot2^{*}_{\mu }}|G^{-1}(v^+(x))|^{\al\cdot2^{*}_{\mu}-2}G^{-1}(v^+(x))w_n}{|x-y|^ \mu \,g(G^{-1}(v^+(x)))}dxdy\\ &\quad+\iint\limits_{\RR^{2N}}\frac{|G^{-1}(w_n^+(y))|^{\al\cdot2^{*}_{\mu }}|G^{-1}(w_n^+(x))|^{\al\cdot2^{*}_{\mu}-2}G^{-1}(w_n^+(x))w_n}{|x-y|^ \mu \,g(G^{-1}(w_n^+(x)))}dxdy + o(1)\quad \text{as}\; n \to \infty.
\end{aligned}
\end{equation}
Thus by substituting \eqref{Deq5.4} and \eqref{Deq5.7} in \eqref{Deq5.6}, we see as $n \to \infty$
  \begin{equation}\label{Deq5.8}
  	\begin{aligned}
  		 \int\limits_{\RR^N} |\nabla w_n|^2 =  \iint\limits_{\RR^{2N}}\frac{|G^{-1}(w_n^+(y))|^{\al\cdot2^{*}_{\mu }}|G^{-1}(w_n^+(x))|^{\al\cdot2^{*}_{\mu}-2}G^{-1}(w_n^+(x))w_n}{|x-y|^ \mu \,g(G^{-1}(w_n^+(x)))}dxdy = k + o(1).
  	\end{aligned}
  \end{equation}
Further by the Remark \ref{Drem1.1} and Lemma \ref{Dlem3.1}(6), we deduce that
\begin{equation}\label{Deq5.9}
	\begin{aligned}
		\frac{|G^{-1}(w_n^+(x))|^{\al\cdot2^{*}_{\mu}-2}G^{-1}(w_n^+(x))w_n}{g(G^{-1}(w_n^+(x)))} \leq \frac{|G^{-1}(w_n^+(x))|^{\al\cdot(2^{*}_{\mu}-1)}w_n}{\ba} \leq \frac{\al^{2^*_{\mu-1}}}{\ba^{2^*_\mu}}|w_n|^{2^*_\mu}.
	\end{aligned}
\end{equation}
Hence by using \eqref{Deq5.8}, \eqref{Deq5.9} and the Sobolev best constant, we get
\begin{equation}\label{Deq5.10o}
\begin{aligned}
	k +o(1) =&	\iint\limits_{\RR^{2N}}\frac{|G^{-1}(w_n^+(y))|^{\al\cdot2^{*}_{\mu }}|G^{-1}(w_n^+(x))|^{\al\cdot2^{*}_{\mu}-2}G^{-1}(w_n^+(x))w_n}{|x-y|^ \mu \,g(G^{-1}(w_n^+(x)))}dxdy \leq \frac{\al^{2\cdot2^*_{\mu-1}}}{\ba^{2\cdot2^*_\mu}}\|w_n\|_0^{2\cdot2^*_\mu}\\ &\leq \frac{\al^{2\cdot2^*_{\mu-1}}}{\ba^{2\cdot2^*_\mu}(S^H)^{2^*_\mu}}\left[ \int\limits_{\RR^N}|\nabla w_n|^2\right] ^{2^*_\mu} = \frac{\al^{2\cdot2^*_{\mu-1}}}{\ba^{2\cdot2^*_\mu}(S^H)^{2^*_\mu}} k^{2^*_\mu},
\end{aligned}
\end{equation}
which implies
\begin{equation}\label{Deq5.10}
	k \geq \frac{1}{\al}\left( \frac{\ba^2 S^H}{\al}\right)^{\frac{2^*_{\mu}}{2^*_{\mu}-1}}.
\end{equation}

Also by Lemma \ref{Dlem3.1}(4), 
we conclude that
\begin{equation}\label{Deq5.11}
	\iint\limits_{\RR^{2N}}\frac{|G^{-1}(w_n^+(y))|^{\al\cdot2^{*}_{\mu }}|G^{-1}(w_n^+(x))|^{\al\cdot2^{*}_{\mu}-2}G^{-1}(w_n^+(x))w_n}{|x-y|^ \mu \,g(G^{-1}(w_n^+(x)))}dxdy \geq \frac{1}{\al}\|G^{-1}(w_n^+)\|_{0,\al},
\end{equation}
Putting together \eqref{Deq5.10} and \eqref{Deq5.11}, in \eqref{Deq5.5}, we get
\begin{equation*}
	c\geq J(v) + \left[ \frac{1}{2}- \frac{1}{2\cdot2^*_{\mu}}\right] k \geq \frac{1}{\al}\left[ \frac{1}{2}- \frac{1}{2\cdot2^*_{\mu}}\right]\left( \frac{\ba^2 S^H}{\al}\right)^{\frac{2^*_{\mu}}{2^*_{\mu}-1}} = c^\infty_{*}
\end{equation*}
and from \eqref{Deq5.1n}, we get a contradiction. Hence the case $\|w_n\|_{L^2} \to 0$ as $n \to \infty$, is impossible.\\
\textbf{Step 2:} When $\|w_n\|_{L^2} \to b >0$ as $n \to \infty$.\\
Using the concentration-compactness principle, we prove that this case is also not possible for $ c \in \left(0, J^\infty \right) $.\\
We proceed using Lemma \ref{Dlem3.5}, with $\rho_n = |w_n|^2$, thus there are two possibilities:\\
\textbf{Case (i):} Vanish occurs\\
\begin{equation*}
	\ds\lim\limits_{n \to +\infty}\sup\limits_{y \in \RR^N}\int\limits_{B_R(y)}|w_n|^2 (x)\,dx = 0,\; \text{for all}\; 0 <R< \infty.
\end{equation*}
Thus by Lemma \ref{Dlem3.6}, we conclude that $w_n \to 0 $ in $L^q(\RR^N)$ for $q \in (2, 2^*)$. Now arguing as in Step 1, we get as $n \to \infty$
\begin{equation}\label{Deq5.12}
	\begin{aligned}
		c +o(1) &= J(v_n)\\
		&= \frac{1}{2}\int\limits_{\RR^N} |\nabla w_n|^2 + a(x)w_n^2 dx -\frac{1}{2\al\cdot{2^*_{\mu } }}\|G^{-1}(w_n^+)\|_{0,\al}^{2\al\cdot2^{*}_{\mu }} + J(v),\\
	\end{aligned}
\end{equation}
and
\begin{equation}\label{Deq5.13}
	\begin{aligned}
		o(1) &= \left\langle J'(v_n) - J'(v), v_n -v\right\rangle \\
		& = \int\limits_{\RR^N} |\nabla w_n|^2 + a(x)w_n^2 -\iint\limits_{\RR^{2N}}\frac{|G^{-1}(w_n^+(y))|^{\al\cdot2^{*}_{\mu }}|G^{-1}(w_n^+(x))|^{\al\cdot2^{*}_{\mu}-2}G^{-1}(w_n^+(x))w_n}{|x-y|^ \mu \,g(G^{-1}(w_n^+(x)))}dxdy.
	\end{aligned}
\end{equation}
By \eqref{Deq5.13}, as $n \to \infty$ we get 
\begin{equation*}
	\begin{aligned}
		\int\limits_{\RR^N} |\nabla w_n|^2 + a(x)w_n^2 dx  &= \iint\limits_{\RR^{2N}}\frac{|G^{-1}(w_n^+(y))|^{\al\cdot2^{*}_{\mu }}|G^{-1}(w_n^+(x))|^{\al\cdot2^{*}_{\mu}-2}G^{-1}(w_n^+(x))w_n}{|x-y|^ \mu \,g(G^{-1}(w_n^+(x)))}dxdy\\
		& = k + o(1).
	\end{aligned}
\end{equation*}
Again by similar argument as in \eqref{Deq5.10o}, we get
\begin{equation*}
		k +o(1) \leq \frac{\al^{2\cdot2^*_{\mu-1}}}{\ba^{2\cdot2^*_\mu}(S^H)^{2^*_\mu}}\left[ \int\limits_{\RR^N}|\nabla w_n|^2 + a(x)w_n^2\,dx\right] ^{2^*_\mu} = \frac{\al^{2\cdot2^*_{\mu-1}}}{\ba^{2\cdot2^*_\mu}(S^H)^{2^*_\mu}} k^{2^*_\mu},
\end{equation*}
which implies
\begin{equation*}
	k \geq \frac{1}{\al}\left( \frac{\ba^2 S^H}{\al}\right)^{\frac{2^*_{\mu}}{2^*_{\mu}-1}}.
\end{equation*}
Thus from \eqref{Deq5.12}, we get
\begin{equation*}
	c\geq J(v) + \left[ \frac{1}{2}- \frac{1}{2\cdot2^*_{\mu}}\right] k \geq c^\infty_{*},
\end{equation*}
which gives us a contradiction.\\
\textbf{Case (ii):} Non-vanish occurs\\
There exist $\al > 0$, $0 <R< \infty$ and $\{y_n\} \subset \RR^N$ such that \begin{equation}\label{Deq5.15}
	\liminf\limits_{n \to +\infty}\int\limits_{B_R(y_n)}|w_n|^2 (x)dx \geq \mathcal{B} > 0.
\end{equation}
Note that $|y_n| \to\infty$ as $n \to \infty$, or else $\{w_n\}$ is tight, and thus $\|w_n\|_{L^2} \to 0$ as $n \to \infty$,
which contradicts $\|w_n\|_{L^2} \to b >0$ as $n \to \infty$. Set $\overline{w}_n(x) = w_n(x + y_n)$. Therefore there exists $w_0 \in H^1(\RR^N)$ such that up to a sub-sequence
\begin{equation*}
	\overline{w}_n \rightharpoonup w_0 \;\text{weakly in}\; H^1(\RR^N),\;\overline{w}_n \to w_0\;\text{in}\; L^s_{loc}(\RR^N)\;\text{for}\; s\in [1, 2^*),\; \overline{w}_n \to w_0\;\text{a.e. in}\;\RR^N.
\end{equation*}
Let $\phi_n(x) = \phi (x-y_n)$, for any $\phi \in C^\infty_0(\RR^N)$, we claim that:
\begin{enumerate}
	\item [(1)] $\left\langle (J^\infty)'(\overline{w}_n), \phi\right\rangle = \left\langle J'(w_n), \phi_n\right\rangle + o(1), $ as $n \to \infty$
	 \item [(2)] $\{\overline{w}_n\}$ is a $(PS)$ of $J^\infty$. 
\end{enumerate}
For proving (1), we need to show the following
\begin{enumerate}
	\item [(a)] $\int\limits_{\RR^N}a(x)w_n\phi_n\,dx = \int\limits_{\RR^N}\overline{w}_n\phi\,dx +o(1) $, which holds true.
	 As $w_n \rightharpoonup 0$ weakly in $H^1(\RR^N )$ and $\lim\limits_{n \to \infty } a(x + y_n) = 1$, we see that for any $\phi \in C^\infty_0(\RR^N)$
	 \begin{equation*}
	 	\begin{aligned}
	 		\int\limits_{\RR^N}a(x)w_n\phi_n\,dx 
	 		&= \int\limits_{\RR^N}[a(x+y_n)-1]\overline{w}_n\phi\,dx + \int\limits_{\RR^N}\overline{w}_n\phi\,dx\\
	 		&= \int\limits_{\RR^N}\overline{w}_n\phi\,dx +o(1), \;\;\text{as}\; n \to \infty,
	 	\end{aligned}
\end{equation*}
where the last inequality is obtained using the H\"older's inequality, boundedness of sequence $\{w_n\}$ and Lebesgue dominated convergence theorem.
\item [(b)] $\int\limits_{\RR^N}f(x, w_n)\phi_n\,dx = \int\limits_{\RR^N}\overline{f}(\overline{w}_n)\phi\,dx +o(1) \;\;\text{as}\; n \to \infty$,\\
To prove this, we use Lemma \ref{Dlem3.2}(5) and Lebesgue dominated convergence theorem,
\begin{equation*}
	\begin{aligned}
		\int\limits_{\RR^N}f(x, w_n)\phi_n\,dx 
		&= \int\limits_{\RR^N}[f(x+y_n,\overline{w}_n) - \overline{f}(\overline{w}_n)]\phi\,dx + \int\limits_{\RR^N}\overline{f}(\overline{w}_n)\phi\,dx\\
		&= \int\limits_{\RR^N}\overline{f}(\overline{w}_n)\phi\,dx +o(1).
	\end{aligned}
\end{equation*}
\end{enumerate}
Thus from (a) and (b), we get the desired claim. Hence $\{\overline{w}_n\}$ is a $(PS)$ of $J^\infty$, and $w_0$ is a weak solution of $(P^\infty)$.
Next we claim that $w_0\neq 0$.\\
By the assumption in Step (2), we have $\|w_n\|_{L^2} \to b >0$ as $n \to \infty$, thus there exists a sequence $\{y_n\}$ satisfying \eqref{Deq5.15} and
there exists $R > 0$, such that
\begin{equation*}
	\int\limits_{B_R(y_n)}|w_n|^2 (x)dx = \tau + o(1), \; \text{as}\; n \to \infty,
\end{equation*}
where $\tau \in (0, b]$ is a constant. Now if we assume $w_0 = 0$, then
\begin{equation*}
	\int\limits_{B_R(y_n)}|w_n|^2 (x)dx = \int\limits_{B_R(0)}|\overline{w}_n|^2 (x)dx = o(1), \; \text{as}\; n \to \infty,
\end{equation*}
which contradicts the fact that $\tau >0$. Hence $w_0 \neq 0$.\\
We claim that $J(w_n) = J^\infty(\overline{w}_n) +o(1)$, as $n \to \infty$.
For proving the claim, we need to show the following
\begin{enumerate}
	\item [(a)] $\int\limits_{\RR^N}a(x)w_n^2\,dx = \int\limits_{\RR^N}|\overline{w}_n|^2\,dx +o(1) $, which holds true, as $n \to \infty$.\\
	Since $\lim\limits_{n \to \infty } a(x + y_n) = 1$, we have for every $\e>0$, there exists $R > 0$ such that if $|x+ y_n| > R$, $|a(x+y_n) - 1| < \e$.
	\begin{equation*}
		\begin{aligned}
			\int\limits_{\RR^N}a(x)w_n^2\,dx 
			&= \int\limits_{|x+y_n| \leq R}[a(x+y_n)-1]\overline{w}_n^2\,dx + \int\limits_{|x+y_n| > R}[a(x+y_n)-1]\overline{w}_n^2\,dx+ \int\limits_{\RR^N}|\overline{w}_n|^2\\
			&= \int\limits_{\RR^N}|\overline{w}_n|^2\,dx +o(1),
		\end{aligned}
	\end{equation*}
\item [(b)] $\int\limits_{\RR^N}F(x, w_n)\,dx = \int\limits_{\RR^N}\overline{F}(\overline{w}_n)\,dx +o(1)$, as $n \to \infty$.\\
To prove this, we use Lemma \ref{Dlem3.2}(5) and Lebesgue dominated convergence theorem,
\begin{equation*}
	\begin{aligned}
		 \left| \int\limits_{\RR^N}F(x+y_n,\overline{w}_n)-\overline{F}(\overline{w}_n) \,dx\right| 
		&= \left| \int\limits_{\RR^N}\int\limits_{0}^{\overline{w}_n}f(x+y_n,t)-\overline{f}(t) \,dt dx\right|\\
		&\leq \left| \int\limits_{\RR^N}\left( f(x+y_n,\te\overline{w}_n)-\overline{f}(\te\overline{w}_n)\right) \overline{w}_n \, dx\right|\\
		&=o(1),
	\end{aligned}
\end{equation*}
as $n \to \infty$, where $\te \in [0,1]$.
\end{enumerate}
Hence from (a) and (b), we prove the desired claim, $J(w_n) = J^\infty(\overline{w}_n) +o(1)$, as $n \to \infty$.\\
Next we denote $z_n = \overline{w}_n -w_0$. Thus by Brezis-Lieb lemma, we claim that, as $n \to \infty$
\begin{equation*}
	J^\infty(\overline{w}_n) = J^\infty(z_n) + J^\infty(w_0) + o(1)
\end{equation*}
and 
\begin{equation*}
	(J^\infty)'(\overline{w}_n) = (J^\infty)'(z_n) + (J^\infty)'(w_0) + o(1),
\end{equation*}
which implies $(J^\infty)'(\overline{z}_n) = o(1)$, as $n \to \infty$. Consequently we can prove that $J^\infty(z_n) \geq 0$. Finally concluding the above results we get, as $n \to \infty$
\begin{equation*}
	o(1) + c = J(v_n) = J(v) + J(w_n) = J(v) + J^\infty(\overline{w}_n) = J(v) + J^\infty(\overline{z}_n) + J^\infty(w_0) \geq J^\infty(w_0) \geq J^\infty,
\end{equation*}
which gives us a contradiction. Hence the non-vanishing also cannot occur. Hence the case $\|w_n\|_{L^2} \to b >0$ as $n \to \infty$, is also impossible. Therefore, we conclude that $\|w_n\| \to 0$, as $n \to \infty$. \qed
\end{lem}

\begin{lem}\label{Dlem5.2}
 Under the assumptions $(g_0)-(g_1)$, $(h_0)-(h_3)$ and $(a_0)$, we claim that
$$\sup\limits_{t \geq 0} J(tw_R)< J^\infty,$$
for $R$ sufficiently large, where $w_R = w(x - R\nu)$, $w$ and $\nu$ are defined in Lemma \ref{Dlem4.5}.
\proof Since $J(tw_R) \to -\infty$ as $t \to \infty$ uniformly in $R \geq  1$, there exists $t_1 > 0$ such that
\begin{equation*}
	\sup\limits_{t \geq 0}J(tw_R) = \sup\limits_{0 \leq t \leq t_1}J(tw_R).
\end{equation*}
It suffices to show, the result holds true for $0 \leq t \leq t_1$, for R large enough.
From the assumption $(a_0)$, it follows that there exist $a_{*}$, $\tau >0$ such that
\begin{equation*}
	1 - a(x) \geq a_{*} > 0 \;\;\text{for all}\;\; |x| \leq \tau,
\end{equation*}
where $\tau \in(0, 1)$.\\
Then, from \eqref{Deq4.27}, we see that for some positive constant $C$
\begin{equation}\label{Deq5.16}
\begin{aligned}
		\int\limits_{\RR^N}(1 - a(x))w_R^2\,dx &\geq \int\limits_{|x| \leq \tau}a_{*}a_1^2(|x-R\nu|+1)^{-(N-1)}e^{-2|x-R\nu|}\,dx\\
&\geq a_{*}a_1^2(R+2)^{-(N-1)}e^{-2R}\int\limits_{|x| \leq \tau}\,dx\\
&\geq CR^{-(N-1)}e^{-2R}.
\end{aligned}
\end{equation}
Further from Lemma \ref{Dlem3.2}(5) and Lemma \ref{Dlem4.5}, we deduce that
\begin{equation}\label{Deq5.17}
	\begin{aligned}
		 \int\limits_{\RR^N}\overline{F}(tw_R)-F(x,tw_R) \,dx
		&= \int\limits_{\RR^N}\int\limits_{0}^{tw_R}\overline{f}(s)- f(x,s)\,ds dx\\ 
		&\leq \int\limits_{\RR^N}\int\limits_{0}^{tw_R}2Ce^{-\zeta|x|}(\e s + C_\e s^p)\,ds dx\\
		&\leq C\left( \e t_1^2\int\limits_{\RR^N}e^{-\zeta|x|}w_R^2\,dx + C_\e t_1^{p+1} e^{-\zeta|x|}w_R^{p+1}\,dx\right) \\
		& \leq C \left( \e R^{-(N-1)}e^{-2R} + C_\e e^{-\min\{\zeta, p+1\}R}\right).
	\end{aligned}
\end{equation}
Employing \eqref{Deq5.16} and \eqref{Deq5.17}, for $\zeta > 2$, we see that for $R$ large enough,
\begin{equation*}
	\begin{aligned}
		J(tw_R) &\leq J^\infty(tw_R) - \frac{t^2}{2}\int\limits_{\RR^N}(1 - a(x))w_R^2\,dx + \int\limits_{\RR^N}\overline{F}(tw_R)-F(x,tw_R) \,dx\\
		&\leq  J^\infty -CR^{-(N-1)}e^{-2R} +  \e C R^{-(N-1)}e^{-2R} + C_\e e^{-\min\{\zeta, p+1\}R}\\
		&< J^\infty.
	\end{aligned}
\end{equation*}
\qed
\end{lem}
\textbf{Proof of Theorem 1.1:} From Lemma \ref{Dlem3.7}, there exists a $(PS)_{c}$ sequence $\{v_n\}$ in $H^1(\RR^N )$ at the level $c$, as defined in \eqref{Deq3.3n}. 
Consequently, by Lemma \ref{Dlem5.1} and Lemma \ref{Dlem5.2} we get
$0 < \mathfrak{a}  < c \leq \sup\limits_{t \geq 0}J(tw_R) < J^{\infty} = \min \{J^\infty, c^\infty_{*} \}$, for large $R$.
Moreover, summing up all the observations we see the mountain-pass Lemma holds true and it follows that $c$ is the critical value of $J$. Hence, there exists a non-trivial solution say $v \in H^1(\RR^N)$ such that $J(v) = c \geq \mathfrak{a} > 0$. Furthermore, by standard regularity argument, it follows that $v \in C^2(\RR^N)$ and by Lemma \ref{Dlem3.1}(2) and the assumption on $h$, we have
\begin{equation*}
	\begin{aligned}
		-&\De v + a(x)v\\
		  &= a(x)\left( v-\frac{G^{-1}(v)}{g(G^{-1}(v))}\right)  +\frac{h(x,G^{-1}(v)) }{g(G^{-1}(v))}+ \left( I_{\vartheta}*|G^{-1}(v^+)|^{\al\cdot2^*_\mu}\right) \frac{|G^{-1}(v^+)|^{\al\cdot2^*_\mu-2}G^{-1}(v^+)}{g(G^{-1}(v^+))}\\
		  &\geq 0
\end{aligned}\end{equation*}
Thus from strong maximum principle, it follows that $v$ is a positive solution of problem \eqref{Deq1.6}
\qed

\noindent
\textbf{Acknowledgement.}
The research of first author is  supported by a grant of UGC (India)  with  JRF grant number: June 18-414344.

\end{document}